\documentclass{siamltex}
\usepackage{amsfonts,latexsym}
\usepackage{epsfig}
\usepackage{amsmath}
\usepackage{graphicx}
\usepackage{amsmath}
\usepackage{color}
\usepackage{ulem}
\usepackage{algorithm,algorithmic}
\newtheorem{remark}{Remark}[section]

\title{A RATIONAL KRYLOV METHODS FOR LARGE-SCALE LINEAR
	MULTIDIMENSIONAL DYNAMICAL SYSTEMS}
\author{
H. Barkouki \thanks{Ecole Nationale Supérieure d'Arts et Métiers, Meknès, Morocco. Email: h.barkouki@umi.ac.ma} \and K. Jbilou \thanks{Université du Littoral, C\^ote d’Opale, batiment H. Poincarré,
50 rue F. Buisson, F-62280 Calais Cedex, France. University UM6P, Benguerir Morocco.  Email: Khalide.Jbilou@univ-littoral.fr} }
\date{}

\begin{document}

\maketitle

\begin {abstract}
In this paper, we investigate the use  of multilinear algebra for reducing the order of multidimensional linear time-invariant (MLTI) systems. Our main tools are tensor rational Krylov subspace methods, which enable us to approximate the system's solution within a low-dimensional subspace. We introduce the tensor rational block Arnoldi and tensor rational block Lanczos algorithms. By utilizing these methods, we develop a model reduction approach based on projection techniques. Additionally, we demonstrate how these approaches can be applied to large-scale Lyapunov tensor equations, which are critical for the balanced truncation method, a well-known technique for order reduction. An adaptive method for choosing the interpolation points is also introduced. Finally, some numerical experiments are reported to show the effectiveness of the proposed adaptive approaches.
\end {abstract}

\begin{keywords} Einstein product, Model recuction, Multi-linear dynamical systems, Tensor Rational  Krylov subspaces, Lyapunov tensor equations.
\end{keywords}

\section{Introduction}

In this work, we develop methods based on projection onto block rational Krylov subspaces, for reducing the order of multidimensional
linear time invariant (MLTI) systems.

Consider the following continuous time MLTI system

\begin{equation}\label{eq1}
	\left\{ 
	\begin{array}{c c c}
		\dot{\mathcal{X}}(t)  & = & \mathcal{A}* \mathcal{X} (t)+\mathcal{B}* \mathcal{U} (t),  \ \ \ \ \mathcal{X}_0=0 \\ 
		\mathcal{Y}(t) & = & \mathcal{C}* \mathcal{X} (t)\\
	\end{array}
	\right. 
	%\end{array}%
	%\]
\end{equation}

where $ \dot{\mathcal{X}}(t) $ is the derivative of the tensor $ \mathcal{X} (t) \in \mathbb{R}^{J_1 \times \ldots \times J_N}, \mathcal{A} \in \mathbb{R}^{J_1 \times \ldots \times J_N \times J_1 \times \ldots \times J_N} $ is a square
tensor, $\mathcal{B} \in \mathbb{R}^{J_1 \times \ldots \times J_N \times K_1 \times \ldots \times K_M} $ and $\mathcal{C} \in \mathbb{R}^{I_1 \times \ldots \times I_M \times J_1 \times \ldots \times J_N}.$ The tensors  $\mathcal{U} (t) \in \mathbb{R}^{K_1 \times \ldots \times K_M}$ and $ \mathcal{Y} (t) \in  \mathbb{R}^{I_1 \times \ldots \times I_M}$ are the control and output
tensors, respectively.

The transfer function associated to the dynamical system (2.2) is given by
\begin{equation}\label{trfct}
	\mathcal{F}(s)=\mathcal{C}* (s \mathcal{I} - \mathcal{A})^{-1}* \mathcal{B},
\end{equation}
see \cite{jbil1, chen1} for more details.

Then, the aim of the model-order reduction problem is to produce a low-order multidimensional system that preserves the important properties of the original system and has
the following form

\begin{equation}\label{eq2}
	\left\{ 
	\begin{array}{c c c}
		\dot{\mathcal{X}}_m (t)  & = & \mathcal{A}_m * \mathcal{X} _m(t)+\mathcal{B}_m * \mathcal{U} (t), \  \\ 
		\mathcal{Y}_m(t) & = & \mathcal{C}_m * \mathcal{X}_m (t)\,\
	\end{array}
	\right. 
	%\end{array}%
	%\]
\end{equation}

%where 

The associated low-order transfer
function is denoted by
\begin{equation}
	\mathcal{F}_m(s)=\mathcal{C}_{m}* (s \mathcal{I}_m - \mathcal{A}_m)^{-1}* \mathcal{B}_m.
\end{equation}

In the context of linear time invariant {\tt LTI} systems, various model reduction methods  have been explored these last
years. Some of them are based on Krylov subspace methods (moment matching)
while others use balanced truncation; see \cite{antoulas, glover, guger1, heyouni17} and the references therein. For the Krylov subspace method, this technique is based on matching the moments of the original transfer function around some selected frequencies to
finding a reduced order model that matches a certain number of moments of the original model around these frequencies.  The standard version of the Krylov algorithms generates reduced-order models that may provide a poor approximation of some frequency dynamics. To address this issue, rational Krylov subspace methods have been developed in recent years \cite{abidi1,gall11,gall12,grim15,grim16,jaimoukha2,jaimoukha1}. However, a significant challenge in these methods is the selection of interpolation points, which need to be carefully chosen to ensure accurate approximations. The same techniques could also be used to construct  reduced order of MLTI systems; see \cite{jbil1}.   In this paper, we are interested in tensor rational  Krylov subspace methods using projections onto special low dimensional  tensor Krylov subspaces via  the Einstein product. This technique  helps to reduce the complexity
of the original system while preserving its essential features, and what we end up with is a lower dimensional reduced MLTI system that has the same structure as the original one. The MLTI systems have been introduced in \cite{chen0}, where the states, outputs and inputs are preserved in a tensor format
and the dynamic is supposed to be described by some multilinear operators. A technique named
tensor unfolding \cite{brezell} which transforms a tensor into a matrix, allows to extend many concepts and notions
in the LTI theory to the tensor case. One of them is the use of the transfer function that describes the
input-output behaviour of the MLTI system, and to quantify the efficiency of the projected reduced
MLTI system. Many properties in the theory of LTI systems such as, reachibility, controllability and
stability, have been generalized in the case of MLTI systems based on the Einstein product \cite{chen0}.

In this paper, we will focus on projection techniques that utilize rational Krylov tensor based subspaces, specifically the tensor rational block Krylov subspace methods using Arnoldi and Lanczos algorithms. This approach will be introduced  by employing
the Einstein product and multilinear algebra to construct a lower-dimensional reduced MLTI system via projection. We will also provide a brief overview of the tensor balanced truncation method. In the second part of the paper, we will examine large-scale tensor equations, such as the continuous-time Lyapunov tensor equation, which arises when applying the tensor balanced truncation model reduction method. We will demonstrate how to derive approximate solutions using tensor-based Krylov subspace methods.

The paper is structured as follows. Section 2 introduces the general notations and provides the key definitions used throughout the document. 
%In Section 3, we present an overview of MLTI systems and describe the tensor classic block Krylov subspaces employed to construct high-fidelity reduced MLTI systems. 
In Section 3, we detail the reduction process via projection using tensor rational block Krylov subspace based on Arnoldi and Lanczos algorithms. We provide some algebraic properties of the proposed
processes and derive explicit formulations of the error between the original and the reduced transfer
functions. An adaptive method for choosing the interpolation points is also introduced in Section 4. In Section 5, we discuss the tensor balanced truncation method. Section 6 outlines the computation of approximate solutions to large-scale discrete-time Lyapunov tensor equations using the tensor rational block Lanczos method. The last section is devoted to some numerical experiments.
%Finally, some
%numerical experiments are reported to show the effectiveness of the proposed methods.

\section{Background and definitions}
In this section, we present some notations and definitions that will be utilized throughout this paper. Unless stated otherwise, we denote tensors with calligraphic letters, while lowercase and uppercase letters are used to represent vectors and matrices, respectively. An N-mode tensor is a multidimensional array represented as $\mathcal{X} \in \mathbb{R}^{J_1 \times \ldots \times J_N}$, with its elements denoted by $x_{j_1,\ldots,j_N}$ with $ 1 \leq j_i \leq J_i, i=1,\ldots,N$. If $N=0$, then $\mathcal{X} \in \mathbb{R}$ is a 0-mode tensor, which we refer to as a scalar; a vector is a 1-mode tensor, and a matrix is a 2-mode tensor. Additional definitions and explanations regarding tensors can be found in \cite{kolda}.

Next we give a definition of the tensor Einstein product, a multidimensional generalization
of matrix product. For further details, refer to \cite{brezell, chen1}.

\begin{definition}
Let  $\mathcal{A} \in \mathbb{R}^{J_1 \times \ldots \times J_N \times K_1 \times \ldots \times K_M}$	and  $\mathcal{B} \in \mathbb{R}^{K_1 \times \ldots \times K_M \times I_1 \times \ldots \times I_L}$ be two tensors. The Einstein product $ A * _{M} B \in \mathbb{R}^{J_1 \times \ldots \times J_N \times I_1 \times \ldots \times I_L}$ is defined by
\begin{equation} \label{Einstprod}
(A \ast_M B) _{j_1 \ldots j_N i_1 \ldots i_L}= \sum _{k_1=1}^{K_1} \sum _{k_2=1}^{K_2} \ldots \sum _{k_M=1}^{K_M} a_{j_1 \ldots j_N k_1 \ldots k_M} b_{k_1 \ldots k_M i_1 \ldots i_L}.
\end{equation}
\end{definition}

The following  special tensors will be considered in this paper: a nonzero tensor 
\( \mathcal{D} = (d_{j_1 \ldots j_N i_1 \ldots i_N}) \in \mathbb{R}^{J_1 \times \ldots \times J_N \times J_1 \times \ldots \times J_N} \) 
is called a diagonal tensor if all the entries \( d_{j_1 \ldots j_N  i_1 \ldots i_N} \) are equal to zero except for the diagonal entries denoted by \( d_{j_1 \ldots j_N j_1 \ldots j_N} \neq 0 \).

If all the diagonal entries are equal to 1 (i.e., \( d_{j_1 \ldots j_N j_1 \ldots j_N} = 1 \)), then 
\( \mathcal{D} \) is referred to as the identity tensor, denoted by \( \mathcal{I} \). Let \( \mathcal{A} \in \mathbb{R}^{J_1 \times \ldots \times J_N \times K_1 \times \ldots \times K_M} \) and 
\( \mathcal{B} \in \mathbb{R}^{K_1 \times \ldots \times K_M \times J_1 \times \ldots \times J_N} \) be two tensors such that
\( b_{k_1 \ldots k_M j_1 \ldots j_N} = a_{j_1 \ldots j_N k_1 \ldots k_M} \). In this cas, 
\( \mathcal{B} \) is called the transpose of \( \mathcal{A} \) and denoted by \( \mathcal{A}^T \). The inverse of a square tensor is defined as follows:

\begin{definition}
A square tensor \( \mathcal{A} \in \mathbb{R}^{J_1 \times \ldots \times J_N \times J_1 \times \ldots \times J_N} \) is invertible if and only if there exists a tensor \( \mathcal{X} \in \mathbb{R}^{J_1 \times \ldots \times J_N \times J_1 \times \ldots \times J_N} \) such that
\[
\mathcal{A} \ast_N \mathcal{X} = \mathcal{X} \ast_N \mathcal{A} = \mathcal{I} \in \mathbb{R}^{J_1 \times \ldots \times J_N \times J_1 \times \ldots \times J_N}.
\]
In that case, \( \mathcal{X} \) is the inverse of \( \mathcal{A} \) and is denoted by \( \mathcal{A}^{-1} \).	
	
\end{definition}

\begin{proposition}
	Consider \( \mathcal{A} \in \mathbb{R}^{J_1 \times \ldots \times J_N \times K_1 \times \ldots \times K_M} \) and \( \mathcal{B} \in \mathbb{R}^{K_1 \times \ldots \times K_M \times I_1 \times \ldots \times I_L} \). Then we have the following relations:
	\begin{enumerate}
		\item \( (\mathcal{A} \ast_M \mathcal{B})^T = \mathcal{B}^T \ast_M \mathcal{A}^T \).
		\item \( \mathcal{I}_M \ast_M \mathcal{B} = \mathcal{B} \) and \( \mathcal{B} \ast_N \mathcal{I}_N = \mathcal{B} \).
	\end{enumerate}
\end{proposition}

The trace, denoted by \( \text{tr}(\cdot) \), of a square-order tensor \( \mathcal{A} \in \mathbb{R}^{J_1 \times \ldots \times J_N \times J_1 \times \ldots \times J_N} \) is given by
\begin{equation}
\text{tr}(\mathcal{A}) = \sum_{j_1 \ldots j_N} a_{j_1 \ldots j_N j_1 \ldots j_N}. 
\end{equation}
The inner product of the two tensors \( \mathcal{X}, \mathcal{Y} \in \mathbb{R}^{J_1 \times \ldots \times J_N \times K_1 \times \ldots \times K_M} \) is defined as
\[
\langle \mathcal{X}, \mathcal{Y} \rangle = \text{tr}(\mathcal{X}^T \ast_N \mathcal{Y}).
\]

where \( \mathcal{X}^T \in \mathbb{R}^{K_1 \times \ldots \times K_M \times J_1 \times \ldots \times J_N} \) is the transpose tensor of \( \mathcal{X} \). If \( \langle \mathcal{X}, \mathcal{Y} \rangle = 0 \), we say that \( \mathcal{X} \) and \( \mathcal{Y} \) are orthogonal. The corresponding norm is the tensor Frobenius norm given by
\begin{equation}\label{def3}
\|\mathcal{X}\| = \sqrt{\text{tr}(\mathcal{X}^T \ast_N \mathcal{X})}. 
\end{equation}

We recall the definition of the eigenvalues of a tensor as outlined in \cite{liang,yuchao}.

\begin{definition}\label{eigenv}
	
	Let \( \mathcal{A} \in \mathbb{R}^{J_1 \times \ldots \times J_N \times J_1 \times \ldots \times J_N} \) be a tensor. The complex scalar \( \lambda \) satisfying
	\[
	\mathcal{A} \ast_N \mathcal{X} = \lambda \mathcal{X},
	\]
	is called an eigenvalue of \( \mathcal{A} \), where \( \mathcal{X} \in \mathbb{R}^{J_1 \times \ldots \times J_N} \) is a nonzero tensor referred to as an eigentensor. The set of all the eigenvalues of \( \mathcal{A} \) is denoted by \( \Lambda(\mathcal{A}) \).
\end{definition}

\subsection{Tensor unfolding}

Tensor unfolding is essential  for tensor computations. It simplifies the process that allows  the extension of many concepts from matrices to tensors. This technique involves rearranging the tensor's entries into either a vector or a matrix format.

\begin{definition}\label{def1}
Consider the following transformation \( \Psi : \mathbb{T}_{J_1 \ldots J_N K_1 \ldots K_M} \rightarrow \mathbb{M}_{|J|, |K|} \) with $\Psi(\mathcal{A}) = A$ defined component wise as
	\[
	\mathcal{A}_{j_1 \ldots j_N k_1 \ldots k_M} \rightarrow A_{\text{ivec}(j, J), \text{ivec}(k, K)}.
	\]
	
	where we refer to \( J = \{ J_1, \ldots, J_N \} \), \( K = \{ K_1, \ldots, K_M \} \) and \( \mathbb{T}_{J_1 \ldots J_N K_1 \ldots K_M} \) is the set of all tensors in \( \mathbb{R}^{J_1 \times \ldots \times J_N \times K_1 \times \ldots \times K_M} \). \( \mathbb{M}_{|J|, |K|} \) is the set of matrices in \( \mathbb{R}^{|J| \times |K|} \) where \( |J| = J_1 \cdots J_N \) and \( |K| = K_1 \cdots K_M \). The index mapping \( \text{ivec}(\cdot, \cdot) \) introduced in \cite{rogers} is given by
	\[
	\text{ivec}(j, J) = j_1 + \sum_{l=2}^{N} (j_l - 1) \prod_{u=1}^{l-1} J_u,
	\]
	\[
	\text{ivec}(k, K) = k_1 + \sum_{v=2}^{M} (k_v - 1) \prod_{s=1}^{v-1} K_s,
	\]
	where \( j := \{ j_1, \ldots, j_N \} \) and \( k := \{ k_1, \ldots, k_M \} \) are two sets of subscripts.
\end{definition}

\subsection{MLTI system theory}
In this section, we will present some theoretical results on MLTI systems that generalize certain aspects of LTI dynamical systems.

We consider the continuous-time MLTI system defined in (\ref{eq1}), and a non-zero initial tensor state $\mathcal{X}_0 = \mathcal{X}(t_0)$. The state solution of this system is given by
\begin{equation}
	\mathcal{X}(t) = e^{\mathcal{A}(t - t_0)} \ast_N \mathcal{X}_0 + \int_{t_0}^{t} e^{\mathcal{A}(t - \tau)} \ast_N \mathcal{B} \ast_M \mathcal{U} (\tau) \, d\tau
\end{equation}

\begin{definition}
	Given the MLTI system (\ref{eq1}). The system is called:
	\begin{enumerate}
		\item Asymptotically stable if \( \|\mathcal{X}_k\| \to 0 \) as \( t \to \infty \).
		\item Stable if \( \|\mathcal{X}_k\| \leq \kappa \|\mathcal{X}_0\| \) for some \( \kappa > 0 \).
	\end{enumerate}
\end{definition}

\begin{proposition}
The MLTI dynamical system (\ref{eq1}) is
\begin{itemize}
	\item Asymptotically stable, if and only if \( \mathcal{A} \) is stable \((\sigma(\mathcal{A}) \subset \mathbb{C}^{-1})\).
	\item Stable, if and only if all eigenvalues of \( \mathcal{A} \) have non-positive real parts, and in addition, all pure imaginary eigenvalues have multiplicity one.
\end{itemize}
\end{proposition}

In the same way as the matrix case, the definitions of the reachability and the observability of a continuous
MLTI system are defined below.

\begin{definition}
The controllability and  observability Gramians associated to the MLTI system (\ref{eq1}) are defined, respectively, as
\[
\mathcal{P} = \int_{0}^{\infty} e^{t \mathcal{A}} \ast_N \mathcal{B} \ast_M \mathcal{B}^T \ast_N e^{t \mathcal{A}^T} \, dt;
\]
and

\[
\mathcal{Q} = \int_{0}^{\infty} e^{t \mathcal{A}^T} \ast_N \mathcal{C}^T \ast_M \mathcal{C} \ast_N e^{t \mathcal{A}} \, dt.
\]

\end{definition}

The two Gramians are the unique solutions of the following coupled Lyapunov
matrix equations
\begin{eqnarray}
\mathcal{A} \ast_N \mathcal{P} + \mathcal{P} \ast_N \mathcal{A}^T + \mathcal{B}\ast_M \mathcal{B}^T & = & \mathcal{O}; \label{contr} \\ 
\mathcal{A}^T \ast_N \mathcal{Q} + \mathcal{Q} \ast_N \mathcal{A} + \mathcal{C}^T \ast_M \mathcal{C} & = & \mathcal{O}. \label{obser}
	\end{eqnarray}

\begin{theorem}
The MLTI system (\ref{eq1}) is controllable if and only if the solution $\mathcal{P}$ of (\ref{contr}) is a weakly symmetric positive-definite square tensor. It is observable if and only if the solution $\mathcal{Q}$ of (\ref{obser}) is a weakly symmetric positive-definite square tensor; see \cite{chen2} for more details about these
notions.
	\end{theorem}

\subsection{Block tensors}
In this subsection, we will briefly outline how tensors can be efficiently stored within a single tensor, similar to the approach used with block matrices. For clarity, we will begin by defining  an even-order paired tensor.

\begin{definition}
(\cite{chen1}) A tensor is called an even-order paired tensor if its order is even, i.e., 
(2N for example), and the elements of the tensor are indicated using a pairwise index 
notation, i.e., \( a_{j_1 i_1 \ldots j_N i_N} \) for \( \mathcal{A} \in \mathbb{R}^{J_1 \times I_1 \times \ldots \times J_N \times I_N} \).
\end{definition}

A definition of an n-mode row block tensor, formed from two even-order paired tensors of the same size, is given as follows:

\begin{definition}
(\cite{chen1}) Consider \( \mathcal{A}, \mathcal{B} \in \mathbb{R}^{J_1 \times I_1 \times \ldots \times J_N \times I_N} \) two even-order paired tensors. The n-mode row block tensor, denoted by $
\left| \mathcal{A} \  \mathcal{B} \right| _n \in \mathbb{R}^{J_1 \times I_1 \times \ldots \times J_n \times 2I_n \times \ldots \times J_N \times I_N}
$

is defined by the following:

\begin{equation}
( \left| \mathcal{A} \  \mathcal{B} \right| _n)_{j_1 i_1 \ldots j_N i_N} =
\begin{cases}
a_{j_1 i_1 \ldots j_N i_N}, & j_k = 1, \ldots, J_k, \, i_k = 1, \ldots, I_k \quad \forall k \\
b_{j_1 i_1 \ldots j_N i_N}, & j_k = 1, \ldots, J_k, \, \forall k, \, i_k = 1, \ldots, I_k \quad \text{for } k \neq n, \\
 \   \   \    \   \text{and}  \  i_k = I_k + 1, \ldots, 2 I_k & \text{for } k = n.
\end{cases}
\end{equation}
\end{definition}

In this study, we focus on 4th-order tensors (i.e., N = 2) without specifically addressing paired ones, meaning we won't use the paired index notation. We have observed that the definition remains applicable for both paired and non-paired 4th-order tensors. To clarify this block tensor notation, we will provide some examples below. Let \( \mathcal{A} \) and \( \mathcal{B} \) be two 4th-order tensors in the space \( \mathbb{R}^{J_1 \times J_2 \times I_1 \times I_2} \). By following the two definitions above, we refer to the 1-mode row block tensor by $ \left| \mathcal{A} \  \mathcal{B} \right|_1 \in \mathbb{R}^{J_1 \times 2J_2 \times I_1 \times I_2}.$ We use the MATLAB colon operation \( : \) to explain how to extract the two tensors \( \mathcal{A} \) and \( \mathcal{B} \). We have:

\[
\left| \mathcal{A} \  \mathcal{B} \right|_1(:, 1 : J_2, :, :) = \mathcal{A} \quad \text{and} \quad   \left| \mathcal{A} \  \mathcal{B} \right|_1(:, J_2 + 1 : 2J_2, :, :) = \mathcal{B}.
\]

The notation of the 2-mode row block tensor $ \left| \mathcal{A} \  \mathcal{B} \right|_2 \in \mathbb{R}^{J_1 \times J_2 \times I_1 \times 2I_2}$ described below is the one most commonly used in this paper. We have:

\[
\left| \mathcal{A} \  \mathcal{B} \right|_2(:, :, :, 1 : I_2) = \mathcal{A} \quad \text{and} \quad \   \left| \mathcal{A} \  \mathcal{B} \right|_2(:, :, :, I_2 + 1 : 2I_2) = \mathcal{B}.
\]

For the cas \( N = 1 \) (i.e., \( \mathcal{A} \) and \( \mathcal{B} \) are now matrices, defined as \( A \) and \( B \)), then the notation $
\left| \mathcal{A} \  \mathcal{B} \right|_1 = \left| \mathcal{A} \  \mathcal{B} \right|_2 = [A, B] \in \mathbb{R}^{J_1 \times 2I_1}$ is now defined by the standard block matrices definition.

In the same manner, we can define the 1 or 2-mode column block tensor for the 4th-order tensors based on the definition proposed by Chen et al. \cite{chen0}. The 1-mode column block tensor in \( \mathbb{R}^{2J_1 \times J_2 \times I_1 \times I_2} \) is denoted by $\begin{array}{|c|}
\mathcal{A}\\
\mathcal{B}\\
\end{array} _  1 $ and defined as follows

\[
\begin{array}{|c|}
\mathcal{A}\\
\mathcal{B}\\
\end{array} _  1 (1 : J_1, :, :, :) = \mathcal{A} \quad \text{and} \quad
\begin{array}{|c|}
\mathcal{A}\\
\mathcal{B}\\
\end{array} _  1 (J_1 + 1 : 2J_1, :, :, :) = \mathcal{B}.
\]

While  the 2-mode column block tensor in \( \mathbb{R}^{J_1 \times J_2 \times 2I_1 \times I_2} \), denoted by $\begin{array}{|c|}
\mathcal{A}\\
\mathcal{B}\\
\end{array} _  2$, is  given by

\[
\begin{array}{|c|}
\mathcal{A}\\
\mathcal{B}\\
\end{array} _  2 (:, :, 1 : I_1, :) = \mathcal{A} \quad \text{and} \quad
\begin{array}{|c|}
\mathcal{A}\\
\mathcal{B}\\
\end{array} _  2 (:, :, I_1 + 1 : 2I_1, :) = \mathcal{B}.
\]

Using the Einstein product,  the following proposition can be easily proved  and used during the computational process described in the following sections.

\begin{proposition}
Consider four tensors of the same size \( \mathcal{A} , \mathcal{B} , \mathcal{C} , \mathcal{D}  \in \mathbb{R}^{K_1 \times K_2 \times K_1 \times K_2} \). Then we have

\[
\underbrace{\left| \mathcal{A} \  \mathcal{B} \right|_2}_{ \in \mathbb{R}^{K_1 \times K_2 \times K_1 \times 2K_2} }\quad \ast_2 \quad
\underbrace{\left| \mathcal{C} \  \mathcal{D} \right|_1}_{ \in \mathbb{R}^{K_1 \times 2K_2 \times K_1 \times K_2}} = \mathcal{A} \ast_2 \mathcal{C} + \mathcal{B} \ast_2 \mathcal{D} \in \mathbb{R}^{K_1 \times K_2 \times K_1 \times K_2}.
\]

\[
\underbrace{\left| \mathcal{A} \  \mathcal{B} \right|_1}_{ \in \mathbb{R}^{K_1 \times K_2 \times K_1 \times 2K_2} }\quad \ast_2 \quad
\underbrace{\left| \mathcal{C} \  \mathcal{D} \right|_2}_{ \in \mathbb{R}^{K_1 \times 2K_2 \times K_1 \times K_2}} =\left|  \left| \mathcal{A} \ast_2 \mathcal{C} \ \mathcal{A} \ast_2 \mathcal{D} \right|_2, \left| \mathcal{B} \ast_2 \mathcal{C} \ \mathcal{B} \ast_2 \mathcal{D} \right|_2 \right|_1 \in \mathbb{R}^{K_1 \times 2K_2 \times K_1 \times 2K_2}.
\]
\end{proposition}

Consider \( m \) tensors \( \mathcal {A}_k \in \mathbb{R}^{J_1 \times J_2 \times I_1 \times I_2} \). By sequentially applying the definition of the 2-mode row block tensor, we can construct a single tensor denoted as \( \mathcal{A} \) in the space \( \mathbb{R}^{J_1 \times J_2 \times I_1 \times mI_2} \), which takes the following form:

\begin{equation}
\mathcal{A} =\left| \mathcal{A}_1  \  \mathcal{A}_2  \ldots   \mathcal{A}_m   \right| \in \mathbb{R}^{J_1 \times J_2 \times I_1 \times m I_2}.
\end{equation}

For more  definitions, see  [\cite{chen1}, Definition 4.2]. The process to
create such tensor can be defined recursively by defining

\[
\mathcal{X}_k =\left| \mathcal{X}_{k-1} \  \mathcal{A}_k \right|  \  \   \text{where}  \  \   \mathcal{X}_1=\mathcal{A}_1. 
\]

An explanation using the MATLAB colon operator for the constructed tensor \( \mathcal{A} \) goes as follows:
\[
\mathcal{A}_l = \mathcal{A}(:, :, :, (l - 1)I_2 + 1 : lI_2) \quad \text{for } l = 1, \ldots, m.
\]

In next sections, we will use some tensors of type \( \mathcal{M}_m \in \mathbb{R}^{K_1 \times mK_2 \times K_1 \times mK_2} \) where \( m > 0 \).  For simplicity,  we consider \( m = 3 \) (i.e., \( \mathcal{M}_3 \in \mathbb{R}^{K_1 \times 3K_2 \times K_1 \times 3K_2} \)),  the tensor $\mathcal{M}_3$ is defined as follows
\begin{equation}
\mathcal{M}_3 = \left|  \left|   \left| \mathcal{M}_{1,1} \ \mathcal{M}_{1,2} \right|_2 \  \mathcal{M}_{1,3}\right|_2, \left| \left| \mathcal{M}_{2,1} \ \mathcal{M}_{2,2} \right|_2 \ \mathcal{M}_{2,3} \right|_2, \left|  \left|  \mathcal{M}_{3,1} \ \mathcal{M}_{3,2}\right|_2 \ \mathcal{M}_{3,3} \right|_2 \right|_1, 
\end{equation}

where \( \mathcal{M}_{i,j} \in \mathbb{R}^{K_1 \times K_2 \times K_1 \times K_2} \). To simplify the notation, we will denote $\mathcal{M}_3$  as

\begin{equation}
\mathcal{M}_3 =
\begin{array}{|c c c|}
\mathcal{M}_{1,1} & \mathcal{M}_{1,2} & \mathcal{M}_{1,3} \\

\mathcal{M}_{2,1} & \mathcal{M}_{2,2} & \mathcal{M}_{2,3} \\

\mathcal{M}_{3,1} & \mathcal{M}_{3,2} & \mathcal{M}_{3,3}\\

\end{array}
\in \mathbb{R}^{K_1 \times 3K_2 \times K_1 \times 3K_2}. 
\end{equation}

The tensors \( \mathcal{M}_{i,j} \) can be extracted from the tensor \( \mathcal{M}_m \) using the MATLAB colon operator as follows:

\[
\mathcal{M}_{i,j} = \mathcal{M}_3(:, (i - 1)K_2 + 1 : iK_2, :, (j - 1)K_2 + 1 : jK_2) \quad \text{for } i, j = 1, 2, 3.
\]

For a generalization to \( m \in \mathbb{N} \) (i.e., \( M_m \in \mathbb{R}^{K_1 \times mK_2 \times K_1 \times mK_2} \)), the tensor $\mathcal{M}_n$ is defined as follows:

\begin{equation}
\mathcal{M}_m =
\begin{array}{|c c c c|}
\mathcal{M}_{1,1} & \mathcal{M}_{1,2} & \ldots & \mathcal{M}_{1,m} \\
\mathcal{M}_{2,1} & \mathcal{M}_{2,2} & \ldots & \mathcal{M}_{2,m} \\
\vdots & \vdots & \ddots & \vdots \\
\mathcal{M}_{m,1} & \mathcal{M}_{m,2} & \ldots & \mathcal{M}_{m,m}
\end{array}  
\in \mathbb{R}^{K_1 \times mK_2 \times K_1 \times mK_2}.
\end{equation}

In the same way, the tensor \( \mathcal{M}_{i,j} \) can be described using the MATLAB colon operator as follows:

\[
\mathcal{M}_{i,j} = \mathcal{M}_m(:, (i-1)K_2 + 1 : iK_2, :, (j-1)K_2 + 1 : jK_2) \in \mathbb{R}^{K_1 \times K_2 \times K_1 \times K_2},  \quad \text{for } i, j = 1, 2, \ldots, m.
\]

For simplicity, let us assume that \( m = 2 \) and \( \mathcal{A}, \mathcal{B} \) are two tensors in \( \mathbb{R}^{K_1 \times K_2 \times K_1 \times K_2} \). Using the definitions described above, we can easily prove the following proposition.

\begin{proposition}
\vskip0.5cm
\begin{enumerate}
\item $ \begin{array}{|c c|}
\mathcal{M}_{1,1} & \mathcal{M}_{1,2} \\
\mathcal{M}_{2,1} & \mathcal{M}_{2,2} \\
\end{array}
= \left|
\left|
\mathcal{M}_{1,1}  \  \mathcal{M}_{1,2} \right| _2 \ 
\left|
\mathcal{M}_{2,1}  \  \mathcal{M}_{2,2} \right| _2  \right| _1
=
\left|
\left|
\mathcal{M}_{1,1}  \  \mathcal{M}_{2,1} \right| _1 \
\left|
\mathcal{M}_{1,2}  \  \mathcal{M}_{2,2} \right| _1  \right| _2.$

\vskip0.5cm

\item $ \left|
\mathcal{A}  \  \mathcal{B} \right| _2 \ast_2
\begin{array}{|c c|}
\mathcal{M}_{1,1}  &  \mathcal{M}_{1,2} \\
\mathcal{M}_{2,1}  & \mathcal{M}_{2,2} \\
\end{array}
=
\left| \mathcal{A} \ast_2 M_{1,1} + \mathcal{B} \ast_2 M_{2,1} \quad \mathcal{A} \ast_2 M_{1,2} + \mathcal{B} \ast_2 M_{2,2} \right|_2. $

\vskip0.5cm
\item $
\begin{array}{|c c|}
\mathcal{M}_{1,1}  &  \mathcal{M}_{1,2} \\
\mathcal{M}_{2,1}  & \mathcal{M}_{2,2} \\
\end{array}   \ast_2 \left|
\mathcal{A}  \  \mathcal{B} \right| _1
=
\left|  \mathcal{M}_{1,1} \ast_2 \mathcal{A}  +  \mathcal{M}_{1,2} \ast_2 \mathcal{B}  \quad  \mathcal{M}_{2,1}  \ast_2 \mathcal{A} + \mathcal{M}_{2,2} \ast_2  \mathcal{B} \right|_1. $
\end{enumerate}
\end{proposition}

\vskip0.5cm

\subsection{QR decomposition}
we  first define the notion of upper and lower triangular tensors which will be used later.
\begin{definition}
Let \( \mathcal{U} \) and \( \mathcal{L} \) be two tensors in the space \( \mathbb{R}^{J_1 \times \ldots \times J_N \times I_1 \times \ldots \times I_N} \).

- \( \mathcal{U} \) is an upper triangular tensor if the entries \( u_{j_1,\ldots,j_N,i_1,\ldots,i_N} = 0 \) when \( \text{ivec}(j, J) \geq \text{ivec}(i, I) \).

- \( \mathcal{L} \) is a lower triangular tensor if the entries \( l_{j_1,\ldots,j_N,i_1,\ldots,i_N} = 0 \) when \( \text{ivec}(j, J) \leq \text{ivec}(i, I) \),

where \( \text{ivec}(\cdot, \cdot) \) is the index mapping mentioned in the Definition \ref{def1}.
\end{definition}

Analogously to the QR decomposition in the matrix case \cite{golub}, a similar definition for the decomposition of a tensor \(\mathcal{A} \) in the space \( \mathbb{R}^{J_1 \times \ldots \times J_N \times K_1 \times \ldots \times K_M} \) is defined as follows:

\begin{equation}
\mathcal{A}  = \mathcal{Q}  \ast \mathcal{R}, 
\end{equation}

where \( \mathcal{Q}  \in \mathbb{R}^{J_1 \times \ldots \times J_N \times K_1 \times \ldots \times K_M} \) is orthogonal, i.e., \( Q^T \ast Q = I_K \), and \( \mathcal{R}  \in \mathbb{R}^{K_1 \times \ldots \times K_M \times K_1 \times \ldots \times K_M} \) is an upper triangular tensor.

\subsection{SVD decomposition}
For  a tensor \(\mathcal{A} \in \mathbb{R}^{J_1 \times \ldots \times J_N \times K_1 \times \ldots \times K_M} \).  The Einstein Singular Value Decomposition (SVD) of $\mathcal{A}$ is defined by \cite{liz}:

\begin{equation}
\mathcal{A}  = \mathcal{U}  \ast _N \Sigma  \ast_M \mathcal{V}^T, 
\end{equation}

where \( \mathcal{U}  \in \mathbb{R}^{J_1 \times \ldots \times J_N \times J_1 \times \ldots \times J_N} \)  and \( \mathcal{V}  \in \mathbb{R}^{K_1 \times \ldots \times K_M \times K_1 \times \ldots \times K_M} \) are orthogonals,  and \( \Sigma  \in \mathbb{R}^{J_1 \times \ldots \times J_N \times K_1 \times \ldots \times K_M} \) is a diagonal tensor that
contains the singular values of $\mathcal{A}$ known also as the Hankel singular values of the associated MLTI system \cite{chen2}.

\vskip0.5cm

 Next, we  recall the tensor classic Krylov subspace based on Arnoldi and Lanczos algorithms. For the rest of the paper and for simplicity, we focus on the case where 
$N=M=2$, which means we will consider 4th-order tensors in the subsequent sections. The results can be readily generalized to the cases where $N \geq 2$ and $ M \geq 2$. In the following discussion, unless otherwise specified, we will denote the Einstein product between two 4th-order tensors  $"*_2$" as "*".

\subsection{Tensor classic Krylov subspace} 

\subsubsection{Tensor block Arnoldi algorithm}

Given two tensors  $ \mathcal{A} \in \mathbb{R}^{J_1 \times J_2 \times J_1 \times J_2}$ and $ \mathcal{B} \in \mathbb{R}^{J_1 \times J_2 \times K_1 \times K_2}$. The \( m \)-th tensor block Krylov subspace is defined by

\begin{equation}\label{krylov}
	\mathcal{K}_m(\mathcal{A}, \mathcal{B}) = \text{Range}\{\mathcal{B}, \mathcal{A} \ast \mathcal{B}, \ldots, \mathcal{A}^{\star (m-1)} \ast \mathcal{B}\} \subseteq \mathbb{R}^{J_1 \times J_2 \times K_1 \times K_2},
\end{equation}

where $ \mathcal{A}^{\star k} = \underbrace{\mathcal{A},\ldots,\mathcal{A}}_{k} $. The tensor block Arnoldi algorithm applied to the pair $(\mathcal{A},\mathcal{B})$  generates a sequence of orthonormal  $J_1 \times J_2 \times K_1 \times K_2$ tensors $\lbrace V_{i} \rbrace$  such that 
$$\mathcal{K}_{m}(\mathcal{A},\mathcal{B})={\rm Range} ( V_{1},V_{2},\ldots,V_{m} ).$$

The tensors $V_i$  that are generated by this algorithm  satisfy the orthogonality conditions, i.e.
\begin{equation}\label{equa2.3}
	\left\{ 
	\begin{array}{c c c}
		V_{j}^{T}* V_{i}=\mathcal{O}, \ \textit{if} \ i \neq j,\\
		V_{j}^{T} *V_{i}=\mathcal{I}, \ \textit{if} \ i = j.\\
	\end{array}
	\right. 
\end{equation}
where $\mathcal{O}$ and $\mathcal{I}$ are the zero tensor and the identity tensor, respectively, of $\mathbb{R}^{K_1\times K_2 \times K_1 \times K_2}$. 
Next, we give a version of the tensor block Arnoldi algorithm that was defined in \cite{guide}. The algorithm is summarized as follows.

\begin{algorithm}[h]
	\caption{The tensor block Arnoldi algorithm ({\tt TBA})}
	\label{algArnoldi}
	\begin{enumerate}
		\item \textbf{Input:} $\mathcal{A}\in\mathbb{R}^{J_1 \times J_2 \times J_1 \times J_2},\ \mathcal{B}\in\mathbb{R}^{J_1 \times J_2 \times K_1 \times K_2}$ and a fixed integer m. 
		\item \textbf{Output:}  $\mathcal{V}_{m+1} \in \mathbb{R}^{J_1 \times J_2 \times K_1 \times (m+1) K_2}$ and $\widetilde{\mathcal{H}}_{m} \in \mathbb{R}^{K_1 \times (m+1) K_2 \times K_1 \times m K_2}$.\\
		\item \textbf {Compute QR decomp.} to $\mathcal{B}$ i.e., $\mathcal{B}=V_1 * H_{1,0}$
		%\item \textbf{Initialize:} $\mathcal{V}_{1}=[V_{1}]$.
		\item \textbf{For} $j=1,\ldots,m$
		%\item \ \  \textbf{if} $(j<m)$
		\item \ \ \ \ \textbf{Set} $\mathcal{W}=\mathcal{A}*V_{j}$   
		\item \ \ \ \ \ \ \textbf{For} $i=1,\ldots,j$
		\item \ \ \ \ \ \  \ \ \ $H_{i,j}=V_{i}^{T}*\mathcal{W} $ 
		\item \ \ \ \ \ \  \ \ \ $\mathcal{W}=\mathcal{W}-V_i*H_{i,j}$.
		\item \ \ \ \ \ \ \textbf{EndFor}
		\item \ \ \ \textbf {Compute QR decomp.} to $\mathcal{W}$ i.e., $\mathcal{W}=V_{i+1} * H_{j+1,j}$
		\item  \textbf{endFor}.
	\end{enumerate}
\end{algorithm}

After m steps, we obtain the following decomposition:

\begin{eqnarray}
	\mathcal{A} \ast \mathcal{V}_m & = & \mathcal{V}_{m+1} \ast \widetilde{\mathcal{H}}_m ,\\  \nonumber 
		& = & \mathcal{V}_{m} \ast \mathcal{H}_m + V_{m+1}*H_{m+1,m}*\mathcal{E}_m ^T, \\ 
	\end{eqnarray}
	
	where \( \mathcal{V}_{m+1} \in \mathbb{R}^{J_1 \times J_2 \times K_1 \times (m+1)K_2} \) contains \( V_i \in \mathbb{R}^{J_1 \times J_2 \times K_1 \times K_2} \) for \( i = 1, \ldots, m + 1 \). The tensor \( \widetilde{\mathcal{H}}_m \in \mathbb{R}^{K_1 \times (m+1) K_2 \times K_1 \times mK_2}\) is a block upper Hessenberg tensor  whose nonzeros block
		entries $H_{i,j} \in \mathbb{R}^{K_1 \times  K_2 \times K_1 \times K_2}$ are defined by Algorithm \ref{algArnoldi} and defined as

		\[
		\widetilde{\mathcal{H}}_m =
				\begin{array}{|c c c c|}
						H_{1,1} & H_{1,2} & \cdots & H_{1,m} \\
						H_{2,1} & H_{2,2} & \cdots & H_{2,m} \\
						\vdots & \vdots & \ddots & \vdots \\
						0 & \cdots & H_{m,m-1} & H_{m,m} \\
						0 & \cdots & 0 & H_{m+1,m}
					\end{array}
				\]

			 where the notation $ \vert \ . \ \vert $ is the definition of block tensors given in the previous paragraph. Finally, the tensor $ \mathcal{E}_m$ is obtained from the identity tensor $\mathcal{I}_m \in \mathbb{R}^{K_1 \times m K_2 \times K_1 \times m K_2} $ as $\mathcal{E}_m=\mathcal{I}_m(:,:,:,(m-1)K_2 +1:mK_2)$.

 \subsection{Tensor block Lanczos algorithm}
Let $\mathcal{B}$ and $\mathcal{C}^T$ be  two initial tensors of $\mathbb{R}^{J_1 \times J_2 \times K_1 \times K_2}$,  and consider the following tensor block Krylov subspaces 
$$\mathcal{K}_m(\mathcal{A}, \mathcal{B}) = \text{Range}\{\mathcal{B}, \mathcal{A} \ast \mathcal{B}, \ldots, \mathcal{A}^{\star (m-1)} \ast \mathcal{B}\} \; {\rm  and} \;  \mathcal{K}_m(\mathcal{A}^T, \mathcal{C}^T) = \text{Range}\{\mathcal{C}^T, \mathcal{A}^T \ast \mathcal{C}^T, \ldots, (\mathcal{A}^{T})^{\star (m-1)} \ast \mathcal{C}^T\}.$$ 

The nonsymmetric tensor block Lanczos algorithm applied to the pairs $(\mathcal{A},\mathcal{B})$ and $(\mathcal{A}^{T},\mathcal{C}^T)$ generates two sequences of bi-orthonormal  $J_1 \times J_2 \times K_1 \times K_2$ tensors $\lbrace V_{i} \rbrace$ and $\lbrace W_{j} \rbrace$ such that 
$$\mathcal{K}_{m}(\mathcal{A}, \mathcal{B})={\rm Range} ( V_{1},V_{2},\ldots,V_{m} ).$$
and 
$$\mathcal{K}_{m}(\mathcal{A}^T, \mathcal{C}^T)={\rm Range} ( W_{1},W_{2},\ldots,W_{m} ).$$
The tensors $V_i$ and $W_j$  that are generated by the tensor block Lanczos algorithm  satisfy the biorthogonality conditions, i.e.
\begin{equation}
\left\{ 
\begin{array}{c c c}
W_{j}^{T} *V_{i}=\mathcal{O}, \ \textit{if} \ i \neq j,\\
 W_{j}^{T} * V_{i}=\mathcal{I}, \ \textit{if} \ i = j.\\
\end{array}
\right. 
\end{equation}

Next, we give a stable version of the nonsymmetric tensor block Lanczos process. This algorithm  is analogous to the one defined in \cite{bai1} for matrices and is summarized as follows.\\
\begin{algorithm}[h]
	\caption{The tensor block Lanczos algorithm ({\tt TBL})}
	\label{alglanc}
\begin{enumerate}
\item \textbf{Inputs:} $\mathcal{A}\in\mathbb{R}^{J_1 \times J_2 \times J_1 \times J_2},\ \mathcal{B}, \mathcal{C}^T \in\mathbb{R}^{J_1 \times J_2 \times K_1 \times K_2}$ and m an integer.
\item \textbf{Compute the QR decomposition of $\mathcal{C} * \mathcal{B}$, i.e.,} $\mathcal{C} * \mathcal{B}= \delta * \beta $;\\
 $V_{1}=\mathcal{B}* \beta^{-1}; W_{1}=\mathcal{C}^T* \delta; \widetilde{V}_{2}=\mathcal{A}*V_{1}; \widetilde{W}_{2}=\mathcal{A}^{T}* W_{1};$
\item \textbf{For} $j=1,\ldots,m$ \\
\quad $\alpha_{j}=W_{j}^{T}* \widetilde{V}_{j+1}; \widetilde{V}_{j+1}=\widetilde{V}_{j+1}-V_{j} *\alpha_{j}; \widetilde{W}_{j+1}=\widetilde{W}_{j+1}-W_{j} *\alpha_{j}^{T};$
\item  \textbf{Compute the QR decomposition of $\widetilde{V}_{j+1}$ and $\widetilde{W}_{j+1}$, i.e.,}

\quad $\widetilde{V}_{j+1}=V_{j+1}* \beta_{j+1}; \widetilde{W}_{j+1}=W_{j+1} *\delta_{j+1}^{T};$
\item  \textbf{Compute the singular value decomposition of $W_{j+1}^{T} *V_{j+1}$, i.e.,}\\
\quad $W_{j+1}^{T}*V_{j+1}=\mathcal{U}_{j} * \Sigma_{j} * \mathcal{Z}_{j}^{T}$;\\
\quad $\delta_{j+1}=\delta_{j+1}* \mathcal{U}_{j}* \Sigma^{1/2}_{j}; \beta_{j+1}= \Sigma^{1/2}_{j}* \mathcal{Z}_{j}^{T}* \beta_{j+1}$;\\
\quad $V_{j+1}=V_{j+1}*\mathcal{Z}_{j} *\Sigma^{-1/2}_{j}; W_{j+1}=W_{j+1}*\mathcal{U}_{j} * \Sigma^{-1/2}_{j};$\\
\quad $\widetilde{V}_{j+2}=\mathcal{A}*V_{j+1}-V_{j}*\delta_{j+1}; \widetilde{W}_{j+2}=\mathcal{A}^{T}*W_{j+1}-W_{j}*\beta_{j+1}^{T};$
\item \textbf{end For}.
\end{enumerate}
\end{algorithm}

Setting $\mathcal{V}_{m}=\left|  V_{1}\ V_{2} \ \ldots \ V_{m} \right| $ and $\mathcal{W}_{m}=\left|  W_{1} \ W_{2} \ \ldots \ W_{m} \right| $, we have the following tensor block Lanczos relations
$$\mathcal{A}*\mathcal{V}_{m}=\mathcal{V}_{m} *\mathcal{T}_{m}+V_{m+1}* \beta_{m+1}* \mathcal{E}_{m}^{T},$$
and
$$\mathcal{A}^{T}*\mathcal{W}_{m}=\mathcal{W}_{m}*\mathcal{T}_{m}^{T}+W_{m+1}* \delta_{m+1}^{T}*\mathcal{E}_{m}^{T},$$
where $\mathcal{T}_m$ is the $K_1 \times mK_2 \times K_1 \times mK_2$ block tensor 
 defined by
\[
\mathcal{T}_{m}=
\begin{array}{|c c c c c|}
\alpha_{1} & \delta_{2}&  &  & \\ 
\beta_{2} & \alpha_{2} & . &  &  \\ 
 &  .& . & .&  \\ 
 &  &  .& . &  \delta_{m} \\ 
 &  &  & \beta_{m}& \alpha_{m}\\
\end{array},
\]

whose nonzeros  entries  $\alpha_{j},  \beta_{j}, \delta_{j}, j=1,\ldots,m$ are block tensors of $ \in \mathbb{R}^{K_1 \times  K_2 \times K_1 \times K_2}$ and defined by Algorithm \ref{alglanc}.

\section{Projection based tensor rational block methods}
\subsection{Tensor rational block Arnoldi method}
%\subsubsection{Rational tensor block Arnoldi algorithm}
Let $\mathcal{A}$ and $\mathcal{B}$ be two tensors of appropriate
dimensions, then the tensor rational block Krylov subspace denoted by $\mathcal{K}_m(\mathcal{A}, \mathcal{B}, \Sigma_m)$,
can be defined as follows

\begin{equation}\label{ratkrylov}
	\mathcal{K}_m(\mathcal{A}, \mathcal{B}, \Sigma_m)= {\rm Range} \lbrace (\mathcal{A}-\sigma_{1}I_n)^{-1}* \mathcal{B}, \ldots,\prod_{k=1}^m (\mathcal{A}-\sigma_{k}I_n)^{-1}* \mathcal{B} \rbrace,  \\
\end{equation}

where $$\prod_{k=1}^m (\mathcal{A}-\sigma_{k}I_n)^{-1}= (\mathcal{A}-\sigma_{1}I_n)^{-1}*\ldots *(\mathcal{A}-\sigma_{m}I_n)^{-1},$$
and $ \Sigma_m = \left\lbrace \sigma_1,\ldots,\sigma_m \right\rbrace  $ is the set of interpolation points.

It is worth noting that if we consider the
isomorphism $\Psi$ defined previously to define the matrices $\Psi (\mathcal{A}) = A \in \mathbb{R}^{J_1 J_2 \times J_1 J_2}$ and $\Psi (\mathcal{B}) = B \in \mathbb{R}^{J_1 J_2 \times K_1 K_2}$, then we can generalize all the results obtained in the matrix case to the tensor structure using the Einstein product. Further details are given below. Next, we describe the process for the construction of
the tensor $\mathcal{V}_m$ associated to the  tensor rational block Krylov subspace defined above. The process is
guaranteed via the following rational Arnoldi algorithm. As mentioned earlier, our interest is focused
only on 4th-order tensors, but the following results remain true also for higher order.

The tensor rational block Arnoldi process described below is an analogue to the one proposed for matrices \cite{abidi1, lee}.

\begin{algorithm}[h]
	\caption{The tensor rational block Arnoldi algorithm ({\tt TRBA})}
	\label{alg1}
	\begin{enumerate}
		\item \textbf{Input:} $\mathcal{A}\in\mathbb{R}^{J_1 \times J_2 \times J_1 \times J_2},\ \mathcal{B}\in\mathbb{R}^{J_1 \times J_2 \times K_1 \times K_2}$ and a fixed integer m. 
		\item \textbf{Output:}  $\mathcal{V}_{m+1} \in \mathbb{R}^{J_1 \times J_2 \times K_1 \times (m+1) K_2}$ and $\widetilde{\mathcal{H}}_{m} \in \mathbb{R}^{K_1 \times (m+1) K_2 \times K_1 \times m K_2}$.\\
		\item \textbf {Compute QR decomp.} to $\mathcal{B}$ i.e., $\mathcal{B}=V_1 * H_{1,0}$
		%\item \textbf{Initialize:} $\mathcal{V}_{1}=[V_{1}]$.
		\item \textbf{For} $j=1,\ldots,m$
		%\item \ \  \textbf{if} $(j<m)$
		\item \ \ \ \ \textbf{Set} $\mathcal{W}=(\mathcal{A}-\sigma_{j}\mathcal{I}_n)^{-1}*V_{j}$   
		\item \ \ \ \ \ \ \textbf{For} $i=1,\ldots,j$
		\item \ \ \ \ \ \  \ \ \ $H_{i,j}=V_{i}^{T}*\mathcal{W} $ 
		\item \ \ \ \ \ \  \ \ \ $\mathcal{W}=\mathcal{W}-V_i*H_{i,j}$.
		\item \ \ \ \ \ \ \textbf{EndFor}
		\item \ \ \ \textbf {Compute QR decomp.} to $\mathcal{W}$ i.e., $\mathcal{W}=V_{i+1} * H_{j+1,j}$
		%\item \  \  \textbf{else}
		%\item \ \ \ \ \textbf{Set} $\mathcal{W}=A^{-1}*\mathcal{B}$   
	%	\item \ \ \ \ \ \ \textbf{For} $=i,\ldots,m$
		%\item \ \ \ \ \ \  \ \ \ $H_{i,m}=V_{i}^{T}*\mathcal{W} $ 
		%\item \ \ \ \ \ \  \ \ \ $\mathcal{W}=\mathcal{W}-V_i*H_{i,m}$.
		%\item \ \ \ \ \ \ \textbf{EndFor}
		%\item \ \ \ \textbf {Compute QR decomp.} to $\mathcal{W}$ i.e., $\mathcal{W}=V_{m+1} * H_{m+1,m}$
	%	\item \  \  \textbf{endif}
		\item  \textbf{endFor}.
	\end{enumerate}
\end{algorithm}

%\newpage

In Algorithm \ref{alg1}, we assume that we are not given the sequence of
shifts $\sigma_1, \sigma_2, . . . , \sigma_{m}$ and then we need to include the procedure to automatically generate this sequence during the iterations of the process. This adaptive procedure
will be defined in the next sections. In this Algorithm, step 5 is used to generate the next Arnoldi tensors $V_{k+1}, k=1,\ldots,m$. To ensure that these block tensors  generated in each iteration are orthonormal, 
we compute the QR decomposition of $\Phi (\mathcal{W})$  (step 10), where $\Phi $ is the isomorphism defined in Definition \ref{def1}. The computed tensor $\mathcal{V}_{m+1}$  contains $V_j$ for $j=1,...,m+1$ which form an orthonormal basis of the tensor rational Krylov
subspace defined in \ref{eq1}; i.e

$$\mathcal{V}_{m+1} ^T * \mathcal{V}_{m+1} =\mathcal{I}_{m+1}$$

Where $\mathcal{I}_{m+1}$ is the identity tensor of  $\mathbb{R}^{K_1 \times (m+1) K_2 \times K_1 \times (m+1) K_2}$  . Let $\mathcal{H}_{m}$  be the $\mathbb{R}^{K_1 \times m K_2 \times K_1 \times m K_2}$ block upper Hessenberg tensor whose nonzeros block entries are defined by Algorithm \ref{alg1} (step 7). Let $\widetilde{\mathcal{H}}_{m}$ and $\widetilde{\mathcal{K}}_{m}$ be the $\mathbb{R}^{K_1 \times (m+1) K_2 \times K_1 \times m K_2}$ block upper Hessenberg tensors defined as 
$$ \widetilde{\mathcal{H}}_{m}= \begin{array}{|c|}
		\mathcal{H}_{m} \\
		H_{m+1,m} *\mathcal{E}_{m}^{T} \\
	\end{array} \ \ , \ \ \widetilde{\mathcal{K}}_{m} = \begin{array}{|c |}
		\mathcal{I}_{m}+\mathcal{H}_{m} * \mathbb{D}_{m} \\
		H_{m+1,m} *\sigma_{m}* \mathcal{E}_{m}^{T} \\
	\end{array}, $$

where $\mathcal{D}_{m}$ is the diagonal tensor diag($\sigma_{1} \mathcal{I}_K,\ldots,\sigma_{m} \mathcal{I}_K $)  and $\lbrace \sigma_{1},\ldots,\sigma_{m} \rbrace$ are the set of interpolation points used in Algorithm \ref{alg1}.  After m steps of Algorithm \ref{alg1}, and specifically
for the extra interpolation points at $\sigma_m=0$, we have:

\begin{eqnarray}\label{condition}
A *\mathcal{V}_{m+1}* \widetilde{\mathcal{H}}_{m} & = & \mathcal{V}_{m+1} * \widetilde{\mathcal{K}}_m = \mathcal{V}_{m} * \mathcal{K}_m \nonumber\\
\mathcal{A}_m = \mathcal{V}_m ^T *A *\mathcal{V}_m  & = & [\mathcal{K}_m  - \mathcal{V}_m ^T *A*V_{m+1}*H_{m+1,m}*\mathcal{E}_m^T]*\mathcal{H}_m^{-1}\nonumber\\
\mathcal{A} *\mathcal{V}_{m} & = & \mathcal{V}_{m}* \mathcal{A}_m + \Gamma _{m}, 
\end{eqnarray} 

where $\Gamma_m=(\mathcal{V}_m*\mathcal{V}_m^{T}-\mathcal{I})*\mathcal{A}*V_{m+1}*H_{m+1,m}*\mathcal{E}_{m}^{T}*\mathcal{H}_m^{-1}.$ 
%The tensors $\widetilde{\mathcal{H}}_{m}$ and $\widetilde {\mathcal{K}}_{m}$ are the block upper Hessenberg tensors of $\mathbb{R}^{K_1 \times (m+1) K_2 \times K_1 \times m K_2}$, given as follows
%\begin{equation*}
%	\widetilde{\mathcal{H}}_{m}=[\widetilde{\mathcal{H}}^{(1)},\widetilde{\mathcal{H}}^{(2)},\ldots,\widetilde{\mathcal{H}}^{(m)} ]\;\;{\rm and} \; \;\\
%	\widetilde{\mathcal{K}}_{m}=[\widetilde{\mathcal{K}}^{(1)},\widetilde{\mathcal{K}}^{(2)},\ldots,\widetilde{\mathcal{K}}^{(m)} ],
%\end{equation*}
%where for $k=1,\ldots,m-1$ the $k$-th block columns are given by 
%$$\widetilde{\mathcal{H}}^{(k)}=  \left[ \begin{array}{c}
%	H_{k} \\
%	H_{k+1,k} \\
%	{\bf 0}\\
%\end{array}\right]  \; {\rm and} \; \widetilde{\mathcal{K}}^{(k)} = \left[ \begin{array}{c}
%	\mathcal{E}_{k} + \sigma_{k+1}  H_{k} \\
%	\sigma_{k+1} H_{k+1,k}  \\
%	{\bf 0}\\
%\end{array}\right ]$$
%and for $k=m$ we have
%$$\widetilde{\mathcal{H}}^{(m)}=  \left[ \begin{array}{c}
%	H_{m} -\mathcal{E}_{1} H_{1,0} \\
%	H_{m+1,m} \\
%\end{array}\right] 
%\;\; {\rm and} \;\; \widetilde{\mathcal{K}}^{(m)} =\left[ \begin{array}{c}
%	-\sigma_{1}  \mathcal{E}_{1} H_{1,0} \\
%	0 \\
%\end{array}\right].$$

%The tensor $ \mathcal{E}_m$ is get it from the identity tensor $\mathcal{I}_m \in \mathbb{R}^{K_1 \times m K_2 \times K_1 \times m K_2} $ as $\mathcal{E}_m=\mathcal{I}_m(:,:,:,(m-1)K_2 +1:mK_2)$. The identity tensors $\mathcal{I}_N $ and $\mathcal{I}_K $ are of $\mathbb{R}^{J_1 \times  J_2 \times J_1 \times  J_2} $ and $\mathbb{R}^{K_1 \times  K_2 \times K_1 \times  K_2} $, respectively.

\subsubsection{Rational approximation}

In this section, we consider tensor rational block Arnoldi algorithm described in
the previous section to approximate the associated transfer function (\ref{trfct}) to the MLTI system (\ref{eq1}). We begin by rewriting this transfer function  as:
$$\mathcal{F}_m(s)=\mathcal{C} * \mathcal{X},$$
where $\mathcal{X}$ verifies the following multi-linear system
\begin{equation}
	(s \mathcal{I} - \mathcal{A})*\mathcal{X}=\mathcal{B}.
\end{equation}

In order to find an approximation to $\mathcal{F}(s)$, it remains to find an approximation to the above multilinear system, which can be done by using a projection into tensor rational Krylov subspace defined in (\ref{ratkrylov}).

Let $\mathcal{V}_m$ be the corresponding basis tensor for the tensor rational Krylov subspace. After approximating the full order state $\mathcal{X}$ by $\mathcal{V}_m * \mathcal{X}$, and analogously to the matrix case, by using the Petrove-Galerkin condition technique, we obtain the desired reduced MLTI system (\ref{eq2}), with the following tensorial structure
\begin{equation}\label{eq3}
	\mathcal{A}_m=\mathcal{V}_{m}^{T}*(\mathcal{A}* \mathcal{V}_{m}), \ \ 	\mathcal{B}_m=\mathcal{V}_{m}^{T}*\mathcal{B}, \ \ \mathcal{C}_m=\mathcal{C}*\mathcal{V}_{m}.
\end{equation}

\subsubsection{Error estimation for the transfer function}
 The computation of the exact error in the transfer tensor between the original and the reduced systems
\begin{equation}\label{eq4}
\varepsilon(s)= \mathcal{F}(s) - \mathcal{F}_m (s)
\end{equation}

is important for the measure of accuracy of the resulting reduced-order model.
Unfortunately, the exact error $\varepsilon (s)$ is not available, because the higher dimension of
the original system yields the computation of $\mathcal{F}(s$) to be very difficult. we propose  the following  simplified expression to the error norm $\Vert \varepsilon \Vert $.

\begin{theorem}
	Let  $\mathcal{F}(s$) and  $\mathcal{F}_m(s$) be the two transfer functions associated to the original
	MLTI system (\ref{eq1}) and the reduced one (\ref{eq2}) respectively. Using  the previous results, we have
	\begin{equation}\label{eq5}
		\Vert \mathcal{F}(s) - \mathcal{F}_m(s)   \Vert  \leq \Vert \mathcal{C} * (s\mathcal{I} - \mathcal{A})^{-1} \Vert * \Vert \mathcal{A}*V_{m+1}*H_{m+1,m}*\mathcal{E}_{m}^T *\mathcal{H}_{m}^{-1}*(s\mathcal{I}_m - \mathcal{A}_m)^{-1}*\mathcal{B}_m \Vert.
	\end{equation}
\end{theorem}

\begin{proof}
The error between the initial and the projected transfer functions is given by
\begin{eqnarray*}
	 \mathcal{F}(s) - \mathcal{F}_m(s) & = &  \mathcal{C}* (s \mathcal{I} - \mathcal{A})^{-1}* \mathcal{B} - \mathcal{C}_{m}* (s \mathcal{I}_m - \mathcal{A}_m)^{-1}* \mathcal{B}_m \nonumber\\
	  & = & \mathcal{C}* (s \mathcal{I} - \mathcal{A})^{-1}*[\mathcal{B}-(s \mathcal{I} - \mathcal{A})*\mathcal{V}_m *  (s \mathcal{I}_m - \mathcal{A}_m)^{-1}* \mathcal{B}_m ]
\end{eqnarray*}
Using decomposition (\ref{condition}) we obtain
\begin{eqnarray*}
\mathcal{F}(s) - \mathcal{F}_m(s) & = &	\mathcal{C}* (s \mathcal{I} - \mathcal{A})^{-1}*[\mathcal{B} - \mathcal{V}_m *\mathcal{B}_m +\Gamma_m*(s \mathcal{I}_m - \mathcal{A}_m)^{-1}* \mathcal{B}_m] \nonumber\\
&= &  \mathcal{C}* (s \mathcal{I} - \mathcal{A})^{-1}*\Gamma_m*(s \mathcal{I}_m - \mathcal{A}_m)^{-1}* \mathcal{B}_m
\end{eqnarray*}

where $\Gamma_m=(\mathcal{V}_m*\mathcal{V}_m^{T}-\mathcal{I})*\mathcal{A}*V_{m+1}*H_{m+1,m}*\mathcal{E}_{m}^{T}*\mathcal{H}_m^{-1}$.

We conclude the proof by  using the fact that $\mathcal{V}_m *\mathcal{B}_m=\mathcal{B} $ and $(\mathcal{V}_m*\mathcal{V}_m^{T}-\mathcal{I})$ is an orthogonal projection.
\end{proof}

\subsection{Tensor rational block Lanczos method}

%\subsubsection{Rational tensor block Lanczos algorithm}
Let $\mathcal{A}, \mathcal{B}$ and $\mathcal{C}$ be three tensors of appropriate
dimensions. The tensor rational block Lanczos procedure  is an algorithm for constructing bi-orthonormal tensors of the union of the block tensor
Krylov subspaces   $\mathcal{K}_m(\mathcal{A}, \mathcal{B})$ and $\mathcal{K}_m(\mathcal{A}^T, \mathcal{C}^T)$. The tensor rational block Lanczos process described below is an analogue to the one proposed for matrices \cite{bar1}.

%\begin{equation}\label{ratkrylov}
%	\mathcal{K}_m(\mathcal{A}, \mathcal{B})= {\rm Range} \lbrace (\mathcal{A}-\sigma_{1}I_n)^{-1}* \mathcal{B}, \ldots,\prod_{k=1}^m (\mathcal{A}-\sigma_{k}I_n)^{-1}* \mathcal{B} \rbrace,  \\
%\end{equation}

 Next, we describe the process for the construction of
the tensors $\mathcal{V}_m$ and $\mathcal{W}_m$ associated to the  tensor rational block Krylov subspace defined above. The process can be
guaranteed via the following rational Lanczos algorithm. As mentioned earlier, our interest is focused
only on 4th-order tensors, but the following results remain true also for higher order.

The tensor rational block Lanczos process described below is an analogue to the one proposed for matrices \cite{abidi1, lee}.

%\vskip0.2cm
%\hrule
%\vskip0.1cm
%\noindent \textbf{\bf Algorithm 2.} The rational tensor block Lanczos algorithm \\
%\vskip0.1cm
%\hrule
%\vskip0.1cm
%\label{alg2}

\begin{algorithm}[h]
	\caption{The tensor rational block Lanczos algorithm ({\tt TRBL})   }
	\label{alg2}
\begin{enumerate}
	\item \textbf{Input:} $\mathcal{A}\in\mathbb{R}^{J_1 \times J_2 \times J_1 \times J_2},\ \mathcal{B}, \mathcal{C}^T \in\mathbb{R}^{J_1 \times J_2 \times K_1 \times K_2}$ and a fixed integer m. 
	%A sequence of interpolation points $\lbrace \sigma_{1},\ldots,\sigma_{m} \rbrace$.
	\item \textbf{Output:} two biorthogonal tensors $\mathcal{V}_{m+1}$ and $\mathcal{W}_{m+1}$ of $\mathbb{R}^{J_1 \times J_2 \times K_1 \times (m+1) K_2}$.\\
%	function $[\mathbb{V}_{m},\mathbb{W}_{m}]$ =\textbf{Rational-Block-Lanczos}(A,B,C,$\lbrace \sigma_{1},\ldots,\sigma_{m} \rbrace$)
	\item \textbf{Set} $\mathcal{S}_{0}=(\mathcal{A}-\sigma_{1} \mathcal{I})^{-1}* \mathcal{B}$\  {\rm and} \ $\mathcal{R}_{0}=(\mathcal{A}-\sigma_{1} \mathcal{I})^{-T} *\mathcal{C}^{T}$  
	\item \textbf{Set} $\mathcal{S}_{0}=V_{1} *H_{1,0}$ \  {\rm and}  \  $\mathcal{R}_{0}=W_{1}*G_{1,0}$ such that $W_{1}^{T}*V_{1}=\mathcal{I}_{K}$;
	\item \textbf{Initialize:} $\mathcal{V}_{1}=[V_{1}]$ \  {\rm and}  \  $\mathcal{W}_{1}=[W_{1}]$.
	\item \textbf{For} $k=1,\ldots,m$
	%\item \ \ \ \textbf{if} $(k<m)$
	\item \ \ \ \ \ \ \textbf{if} $\lbrace\sigma_{k+1} = \infty\rbrace$; $\mathcal{S}_{k}=\mathcal{A}*V_{k}$ \   {\rm and}  \ $\mathcal{R}_{k}=\mathcal{A}^{T}*W_{k};$ \textbf{else}
	\item \ \ \ \ \ \  $\mathcal{S}_{k}=(\mathcal{A}-\sigma_{k+1} \mathcal{I})^{-1}*V_{k}$ \   {\rm and}  \ $\mathcal{R}_{k}=(\mathcal{A}-\sigma_{k+1}\mathcal{I})^{-T}*W_{k};$  \textbf{endif}
	\item \ \ \ \ \ \  $H_{k}=\mathcal{W}_{k}^{T}*\mathcal{S}_{k} $ \  {\rm and} \ $G_{k}=\mathcal{V}_{k}^{T}*\mathcal{R}_{k} $;
	\item \ \ \ \ \ \ $\mathcal{S}_{k}=\mathcal{S}_{k}-\mathcal{V}_{k}*H_{k}$ \  {\rm and} \ $\mathcal{R}_{k}=\mathcal{R}_{k}-\mathcal{W}_{k}*G_{k};$
	\item \ \ \ \ \ \ $\mathcal{S}_{k}=V_{k+1}*H_{k+1,k}$ \  {\rm and} \ $\mathcal{R}_{k}=W_{k+1}*G_{k+1,k};$ \ \ \ \  (QR factorization)
	\item \ \ \ \ \ \ $W_{k+1}^{T}*V_{k+1}=\mathcal{P}_{k}*\mathcal{D}_{k}*\mathcal{Q}_{k}^{T}$; \ \  \ \ (Singular Value Decomposition)
	\item \ \ \ \ \ \ $V_{k+1}=V_{k+1}*\mathcal{Q}_{k}*\mathcal{D}_k^{-1/2}$ \  {\rm and} \ $W_{k+1}=W_{k+1}*\mathcal{P}_{k}* \mathcal{D}_k^{-1/2}$;
	\item \ \ \ \ \ \ $H_{k+1,k}=\mathcal{D}_k^{1/2}*\mathcal{Q}_{k}^{T}*H_{k+1,k}$ \  {\rm and} \ $G_{k+1,k}=\mathcal{D}_k^{1/2}*\mathcal{P}_{k}^{T}*G_{k+1,k}$;
	\item \ \ \ \ \ \ $\mathcal{V}_{k+1}=\left| \mathcal{V}_{k} \ V_{k+1} \right|; \  \mathcal{W}_{k+1}=\left| \mathcal{W}_{k} \ W_{k+1} \right |$; 
%	\item \ \ \ \textbf{else}
%	\item \ \ \ \ \ \ \textbf{if} $\lbrace\sigma_{m+1} = \infty\rbrace$; $\mathcal{S}_{m}=\mathcal{A}*\mathcal{B}$ \   {\rm and}  \ $\mathcal{R}_{m}=\mathcal{A}^{T}*\mathcal{C};$ \textbf{else}
%	\item \ \ \ \ \ \ $\mathcal{S}_{m}=\mathcal{A}^{-1}*\mathcal{B}$ \   {\rm and}\ $\mathcal{R}_{m}=\mathcal{A}^{-T}*\mathcal{C}^{T}$; \textbf{endif}
	%	\item \ \ \ \ \ \  $H_{m}=\mathcal{W}_{m}^{T}*\mathcal{S}_{m} $ \  {\rm and} \ $G_{m}=\mathcal{V}_{m}^{T}*\mathcal{R}_{m} $;
%	\item \ \ \ \ \ \ $\mathcal{S}_{m}=\mathcal{S}_{m}-\mathcal{V}_{m}*H_{m}$ \  {\rm and} \ $\mathcal{R}_{m}=\mathcal{R}_{m}-\mathcal{W}_{m}*G_{m};$
%	\item \ \ \ \ \ \ $\mathcal{S}_{m}=V_{m+1}*H_{m+1,m}$ \  {\rm and} \ $\mathcal{R}_{m}=W_{m+1}*G_{m+1,m};$ \ \ \ \  (QR factorization)
%	\item \ \ \ \ \ \ $W_{m+1}^{T}*V_{m+1}=\mathcal{P}_{m}*\mathcal{D}_{m}*\mathcal{Q}_{m}^{T}$; \ \  \ \ (Singular Value Decomposition)
%	\item \ \ \ \ \ \ $V_{m+1}=V_{m+1}*\mathcal{Q}_{m}*\mathcal{D}_m^{-1/2}$ \  {\rm and} \ $W_{m+1}=W_{m+1}*\mathcal{P}_{m}* \mathcal{D}_m^{-1/2}$;
%	\item \ \ \ \ \ \ $H_{m+1,m}=\mathcal{D}_m^{1/2}*\mathcal{Q}_{m}^{T}*H_{m+1,m}$ \  {\rm and} \ $G_{m+1,m}=\mathcal{D}_m^{1/2}*\mathcal{P}_{m}^{T}*G_{m+1,m}$;
	%\item \ \ \ \ \ \ $\mathcal{V}_{m+1}=[\mathcal{V}_{m},V_{m+1}]; \  \mathcal{W}_{m+1}=[\mathcal{W}_{m},W_{m+1}]$;
%	\item \ \ \ \textbf{endif}
	\item  \textbf{endFor}.
\end{enumerate}
\end{algorithm}
%\vskip0.2cm
%\hrule
%\vskip0.2cm

We notice that  the adaptive procedure to automatically generate the shifts $\sigma_1, \sigma_2, \ldots,
\sigma_{m}$  will be defined in the next sections. In Algorithm \ref{alg2}, steps 7-8 are used to generate the next Lanczos tensors $V_{k+1}$ and $W_{k+1}$. To ensure that theses block tensors  generated in each iteration are biorthogonal, 
we apply the QR decomposition to $\Phi (\mathcal{S}_{k})$ and  $\Phi (\mathcal{R}_{k})$ and then we compute the singular value decomposition of $\Phi(W_{k+1}^{T} \ast V_{k+1}) $ (step  11 and step 12).

The tensors $H_{k}$ and $G_{k}$ constructed in step 9 are $K_1 \times kK_2 \times K_1 \times kK_2$ and they are used to construct the block upper Hessenberg tensors $\mathcal{H}_{m}$ and $\mathcal{G}_{m}$, respectively.

The computed tensors $\mathcal{V}_{m+1}$ and $\mathcal{W}_{m+1}$ from Algorithm \ref{alg2} are bi-orthonormal; i.e.,

$$\mathcal{W}_{m+1} ^T * \mathcal{V}_{m+1} =\mathcal{I}_{m+1},$$
 Let $\mathcal{H}_{m}$ and $\mathcal{G}_{m}$   be the $\mathbb{R}^{K_1 \times m K_2 \times K_1 \times m K_2}$ block upper Hessenberg tensors defined as
$$ \mathcal{H}_{m}= \left |  H_1 \  H_2 \ \ldots \  H_m  \right |  \  \  \text{ and} \  \  \mathcal{G}_{m}= \left |  G_1 \ G_2 \ \ldots \ G_m  \right |.$$
%whose nonzeros block entries are defined by Algorithm \ref{alg1}. 

Let $\widetilde{\mathcal{H}}_{m}, \widetilde{\mathcal{G}}_{m}, \widetilde{\mathcal{K}}_{m}$ and $\widetilde{\mathcal{L}}_{m}$ be the $\mathbb{R}^{K_1 \times (m+1) K_2 \times K_1 \times mK_2}$ block upper Hessenberg tensors defined as 
$$ \widetilde{\mathcal{H}}_{m}= \begin{array}{|c |}
		\mathcal{H}_{m} \\
		H_{m+1,m}*\mathcal{E}_{m}^{T} \\
	\end{array} \ \ , \ \ \widetilde{\mathcal{K}}_{m} =\begin{array}{|c |}
		\mathcal{I}_{m}+\mathcal{H}_{m} * \mathbb{D}_{m} \\
		H_{m+1,m} *\sigma_{m} *\mathcal{E}_{m}^{T} \\
	\end{array}, $$ 

and

$$ \widetilde{\mathcal{G}}_{m}= \begin{array}{|c |}
		\mathcal{G}_{m} \\
		G_{m+1,m}*\mathcal{E}_{m}^{T} \\
	\end{array} \ \ , \ \ \widetilde{\mathcal{L}}_{m}=\begin{array}{|c |}
		\mathcal{I}_{m}+\mathcal{G}_{m} * \mathbb{D}_{m} \\
		G_{m+1,m} *\sigma_{m} *\mathcal{E}_{m}^{T} \\
	\end{array}, $$

where $\mathcal{D}_{m}$ is the diagonal tensor diag($\sigma_{1} \mathcal{I}_K,\ldots,\sigma_{m} \mathcal{I}_K $)  and $\lbrace \sigma_{1},\ldots,\sigma_{m} \rbrace$ are the set of interpolation points used in Algorithm \ref{alg2}.  After m steps of Algorithm \ref{alg2}, and specifically
for the extra interpolation points at $\sigma_m=0$, we have:

	\begin{eqnarray}\label{condit2}
		\mathcal{A} *\mathcal{V}_{m+1}* \widetilde{\mathcal{H}}_{m} & = & \mathcal{V}_{m+1} * \widetilde{\mathcal{K}}_m = \mathcal{V}_{m} * \mathcal{K}_m \nonumber\\
		\mathcal{A}^T *\mathcal{W}_{m+1}* \widetilde{\mathcal{G}}_{m} & = & \mathcal{W}_{m+1} * \widetilde{\mathcal{L}}_m = \mathcal{W}_{m} * \mathcal{L}_m \nonumber\\
		\mathcal{T}_m = \mathcal{W}_m ^T *\mathcal{A} *\mathcal{V}_m  & = & [\mathcal{K}_m  - \mathcal{W}_m ^T *A*V_{m+1}*H_{m+1,m}* \mathcal{E}_m^T]*\mathcal{H}_m^{-1} \nonumber\\
		\mathcal{A} *\mathcal{V}_{m} & = & \mathcal{V}_{m}* \mathcal{T}_m + \Gamma _{m, A}, \nonumber\\
		\mathcal{T}_m ^T = \mathcal{V}_m ^T *\mathcal{A}^T *\mathcal{W}_m  & = & [\mathcal{L}_m  - \mathcal{V}_m ^T *A^T*W_{m+1}*G_{m+1,m}*\mathcal{E}_m^T]*\mathcal{G}_m^{-1} \nonumber\\
		\mathcal{A}^T *\mathcal{W}_{m} & = & \mathcal{W}_{m}* \mathcal{T}_m ^T + \Gamma _{m, A^T}, \nonumber\\
			\end{eqnarray} 
	
	where $$\Gamma _{m, A} = (\mathcal{V}_{m} * \mathcal{W}_{m}^T - \mathcal{I}) * \mathcal{A}* V_{m+1} *H_{m+1,m} * \mathcal{E}_m ^T* \mathcal{H}_m^{-1},$$
	
	and
	
	$$\Gamma _{m, A^T} = (\mathcal{W}_{m} * \mathcal{V}_{m}^T - \mathcal{I}) * \mathcal{A}^T* W_{m+1} *G_{m+1,m} * \mathcal{E}_m ^T* \mathcal{G}_m^{-1}.$$

Let  $\mathcal{V}_{m+1}$ and $\mathcal{W}_{m+1}$ be the bi-orthonormal tensors computed using Algorithm \ref{alg2}. After approximating the full order state $\mathcal{X}$ by $\mathcal{V}_m * \mathcal{X}$, and analogously to the matrix case, by using the Petrove-Galerkin condition technique, we obtain the desired reduced MLTI system  with the following tensorial structure
\begin{equation}\label{eq6}
	\mathcal{A}_m=\mathcal{W}_{m}^{T}*(\mathcal{A}* \mathcal{V}_{m}), \ \ 	\mathcal{B}_m=\mathcal{W}_{m}^{T}*\mathcal{B}, \ \ \mathcal{C}_m=\mathcal{C}*\mathcal{V}_{m}.
\end{equation}

\subsubsection{Error estimation for the transfer function}

In this section, we give a simplified expression to the error norm $\Vert \mathcal{F}(s) - \mathcal{F}_m (s) \Vert $.

\begin{theorem}
	Let  $\mathcal{F}(s$) and  $\mathcal{F}_m(s$) be the two transfer functions associated to the original
	MLTI (\ref{eq1}) and the reduced one (\ref{eq2}) respectively. Using  the previous results, we have
	\begin{equation}\label{eq66}
		\Vert \mathcal{F}(s) - \mathcal{F}_m(s)   \Vert  \leq \Vert \mathcal{C} * (s\mathcal{I} - \mathcal{A})^{-1} \Vert * \Vert \mathcal{A}*V_{m+1}*H_{m+1,m}*\mathcal{E}_{m}^T *\mathcal{H}_{m}^{-1}*(s\mathcal{I}_m - \mathcal{A}_m)^{-1}*\mathcal{B}_m \Vert.
	\end{equation}
\end{theorem}

Next, we present a modeling error in terms of two residual tensors. The results provided below have been established in the matrix case \cite{bar1} and are generalized here for the tensor case.

Let
\begin{equation*}
	\left\{ 
	\begin{array}{c c c}
		\mathcal{R}_{B}(s) & = & \mathcal{B}-(s \mathcal{I}_{n}- \mathcal{A})*\mathcal{V}_{m}*\tilde{\mathcal{X}}_{\mathcal{B}}(s),\\ 
		\mathcal{R}_{C}(s) & = & \mathcal{C}^{T}-(s \mathcal{I}_{n}-\mathcal{A})^{T}*\mathcal{W}_{m}*\tilde{\mathcal{X}}_{\mathcal{C}}(s)\\
	\end{array}
	\right. 
	%\end{array}%
	%\]
\end{equation*}
be the tensor residual expressions, where $\tilde{\mathcal{X}}_{\mathcal{B}}(s)$ and $\tilde{\mathcal{X}}_{\mathcal{C}}(s)$ are the solutions of the tensor  equations
\begin{equation*}
	\left\{ 
	\begin{array}{c c c}
		(s\mathcal{I}_{mp}-\mathcal{A}_m) *\tilde{\mathcal{X}}_{\mathcal{B}}(s) & = & \mathcal{B}_m,\\ 
		(s\mathcal{I}_{mp}-\mathcal{A}_m)^{T}* \tilde{\mathcal{X}}_{\mathcal{C}}(s) & = &\mathcal{C}_m^{T},\\
	\end{array}
	\right. 
	%\end{array}%
	%\]
\end{equation*}
and satisfy the Petrov-Galerkin conditions
\begin{equation*}
	\left\{ 
	\begin{array}{c c c}
		\mathcal{R}_{B}(s) & \bot & {\rm Range}({W}_{1} ,\ldots,{W}_{m} )\\ 
		\mathcal{R}_{C}(s) & \bot & {\rm Range} ({V}_{1} ,\ldots,{V}_{m}) ,\\
	\end{array}
	\right. 
	%\end{array}%
	%\]
\end{equation*}
which means that $ \mathcal{W}_{m}^{T}* \mathcal{R}_{\mathcal{B}}(s)=\mathcal{V}_{m}^{T}* \mathcal{R}_{\mathcal{C}}(s)=\mathcal{O}$. In the following theorem, we give an expression of the error $\epsilon(s)$. 
\begin{theorem}
	The error between the frequency responses of the original and reduced-order systems can be expressed as
	\begin{equation}\label{equa4.1.1}
		\epsilon(s)=\mathcal{R}_{\mathcal{C}}^{T}(s)*(s \mathcal{I}_{n}- \mathcal{A})^{-1} *\mathcal{R}_{\mathcal{B}}(s).
	\end{equation}
\end{theorem}

Next, we use the tensor rational Lanczos equations  ( \ref{condit2}) to simplify the tensor residual error expressions. The expressions of the tensor residual $\mathcal{R}_{\mathcal{B}}(s)$ and $\mathcal{R}_{\mathcal{C}}(s)$ could  be written as 
\begin{eqnarray}\label{equa4.8}
	\mathcal{R}_{\mathcal{B}}(s) & = & \mathcal{B}-(s \mathcal{I}_n-\mathcal{A})* \mathcal{V}_{m}* (s \mathcal{I}_{mp}-\mathcal{A}_m)^{-1}* \mathcal{B}_m \nonumber\\
	& = & \underbrace{(\mathcal{V}_{m} * \mathcal{W}_{m}^{T}-I)*\mathcal{A} *V_{m+1}}_{\tilde{\mathcal{B}}}* \underbrace{H_{m+1,m}* \mathcal{E}_{m}^{T}* \mathcal{H}_m^{-1}*(s \mathcal{I}_{mp}- \mathcal{A}_m)^{-1}* \mathcal{B}_m}_{\tilde{\mathcal{R}}_{\mathcal{B}}(s)} 
\end{eqnarray}
\begin{eqnarray}\label{equa4.9}
	\mathcal{R}_{\mathcal{C}}(s) & = & \mathcal{C}^{T}-(s \mathcal{I}-\mathcal{A})^{T} * \mathcal{W}_{m} (s \mathcal{I}_{mp}-\mathcal{A}_m)^{-T} * \mathcal{C}_m^{T} \nonumber\\
	& = & \underbrace{(\mathcal{W}_{m}* \mathcal{V}_{m}^{T}- \mathcal{I})* \mathcal{A}^{T} *W_{m+1}}_{\tilde{\mathcal{C}}^{T}} \underbrace{G_{m+1,m} * \mathcal{E}_{m}^{T} *\mathcal{G}_{m}^{-1}*(s\mathcal{I}_{mp}- \mathcal{A}_m)^{-T} \mathcal{C}_m^{T}}_{\tilde{\mathcal{R}}_{\mathcal{C}}(s)}, 
\end{eqnarray}
where $\tilde{\mathcal{R}}_{\mathcal{B}}(s)$, $\tilde{\mathcal{R}}_{\mathcal{C}}(s)$ are the  terms of the residual errors $\mathcal{R}_{\mathcal{B}}(s)$ and $\mathcal{R}_{\mathcal{C}}(s)$, respectively, depending on the frequencies. The matrices $\tilde{\mathcal{B}}$, $\tilde{\mathcal{C}}^{T}$ are  frequency-independent terms of $\mathcal{R}_{\mathcal{B}}(s)$ and $\mathcal{R}_{\mathcal{C}}(s)$, respectively.
Therefore, the error expression in (\ref{equa4.1.1}) becomes 

\begin{equation}\label{equa4.11}
	\epsilon(s) =  \tilde{\mathcal{R}}_{\mathcal{C}}(s)^{T}* \tilde{\mathcal{C}}*(s\mathcal{I}_{n}-\mathcal{A})^{-1} *\tilde{\mathcal{B}}* \tilde{\mathcal{R}}_{\mathcal{B}}(s) = \tilde{\mathcal{R}}_{\mathcal{C}}(s)^{T} *\tilde{\mathcal{H}}(s)* \tilde{\mathcal{R}}_{\mathcal{B}}(s).
\end{equation}

The transfer function $\tilde{\mathcal{H}}(s)=\tilde{\mathcal{C}}*(s\mathcal{I}_{n}-\mathcal{A})^{-1} * \tilde{\mathcal{B}}$ include terms related to the original system which makes the computation of $\Vert \tilde{\mathcal{R}}_{\mathcal{C}}^{T} *\tilde{\mathcal{H}}*\tilde{\mathcal{R}}_{\mathcal{B}} \Vert $ quite costly. Therefore, instead of using $\tilde{\mathcal{H}}(s)$ we can employ an  approximation of $\tilde{\mathcal{H}}(s)$.
Various possible approximations of the error $\epsilon(s)$ are summarized in Table \ref{tab3}.
\begin{table}[htbp]
	\begin{center}
		\caption{Estimations of the error  $\epsilon(s)$}\label{tab3}
		\begin{tabular}{c|l}
			\hline
			1& $\hat{\epsilon}(s)= \tilde{\mathcal{R}}_{\mathcal{B}}(s)$  \\
			%\hline
			2& $\hat{\epsilon}(s)= \tilde{\mathcal{R}}_{\mathcal{C}}(s)^{T}$ \\
			%\hline
			3 &  $\hat{\epsilon}(s)= \tilde{\mathcal{H}}_m(s)*\tilde{\mathcal{R}}_{\mathcal{B}}(s)$    \\
			%\hline
			4& $\hat{\epsilon}(s)= \tilde{\mathcal{H}}_m(s)$\\
			%\hline
			5 & $\hat{\epsilon}(s)= \tilde{\mathcal{R}}_{C}^{T}(s)*\tilde{\mathcal{H}}_m(s)$ \\
			%\hline
			6 & $\hat{\epsilon}(s)=\tilde{\mathcal{R}}_{C}^{T}(s)*\tilde{\mathcal{H}}_m(s)*\tilde{\mathcal{R}}_{\mathcal{B}}(s)$\\ 
			\hline
		\end{tabular}
	\end{center}
\end{table}

\section{pole section}
In this section, we derive some techniques for adaptive pole selection by employing the representations of the tensor residuals discussed in the previous section.

The goal of these methods is to construct the next interpolation point at each step, based on the idea that shifts should be chosen to minimize the norm of a specific error approximation at every iteration. In this context, an adaptive approach is suggested, that utilizes the following error approximation expression:

\begin{equation*}
	\hat{\epsilon}(s)=\tilde{\mathcal{R}}_{\mathcal{C}}^{T}(s)*\tilde{\mathcal{R}}_{\mathcal{B}}(s).
	%R_{C}^{T}(s)R_{B}(s). 
\end{equation*}
Then the next shift $\sigma_{k+1} \in \mathbb{R}$ is selected as
\begin{equation}\label{equa5.1}
	{ \sigma}_{k+1}=  {\rm arg}  \max_{s \in S} \Vert \tilde{\mathcal{R}}_{\mathcal{C}}^{T}(s)*\tilde{\mathcal{R}}_{\mathcal{B}}(s) \Vert ,
\end{equation}
and if $ \sigma_{k+1}$ is complex, its real part is retained and used as the next interpolation point.\\

For the case of the tensor rational block Arnoldi algorithm, one can choose the next interpolation point as follows:
	\begin{equation}\label{equa5.2}
		{\sigma}_{k+1}=  {\rm arg}  \max_{s \in S} \Vert \tilde{\mathcal{R}}_{\mathcal{B}}(s) \Vert .
	\end{equation}

\section{Tensor Balanced Truncation}
%Tensor Balanced Truncation (TBT) is an advanced model reduction technique specifically tailored for dynamic systems with multiple state variables. This method stands out due to its ability to handle systems whose dynamics are described by multidimensional tensors, a generalization of traditional matrices. By combining principles from model reduction and tensor theory. Tensor Balanced Truncation aims to preserve the essential characteristics of the system while significantly reducing its dimensionality. This approach finds applications across various fields such as controlling complex systems, modeling neural networks, and simulating high-dimensional physical systems.
Tensor Balanced Truncation (TBT) is an advanced technique in numerical linear algebra and model reduction. It extends the principles of Balanced Truncation, a well-established method for reducing the dimensionality of large-scale linear dynamical systems, to tensor computations. This approach finds applications across various fields such as controlling complex systems, modeling neural networks, and simulating high-dimensional physical systems. The core idea behind Balanced Truncation is to approximate a high-dimensional system with a lower-dimensional one while preserving essential dynamical properties. This is achieved by truncating the system based on its controllability and observability Gramians, thus retaining the most significant states. Tensor Balanced Truncation adapts this framework to higher-order tensors, addressing the challenges posed by their complex structure and significant computational demands.

One of the primary motivations for using TBT lies in its potential to handle large-scale tensor data more efficiently. As datasets grow in complexity and size, traditional matrix-based methods can become infeasible, prompting the need for more sophisticated and scalable approaches. TBT offers a promising solution by leveraging the inherent multi-dimensional structure of tensors, thereby enabling more effective data compression and analysis.

 In this section, we present the tensor Balanced Truncation method. The procedure of this method is as
 follows. First, we need to solve two continuous Lyapunov equations
 
 \begin{equation}\label{gramians}
 	\left\{ 
 	\begin{array}{c c c}
 	 \mathcal{A}* \mathcal{P}+\mathcal{P}*\mathcal{A}^T+\mathcal{B}*\mathcal{B}^T=\mathcal{O}\\ 
 	 \mathcal{A}^T*\mathcal{Q}+\mathcal{Q}*\mathcal{A}+\mathcal{C}^T*\mathcal{C}=\mathcal{O},\\
 	\end{array}
 	\right. 
 	%\end{array}%
 	%\]
 \end{equation}
 
where $\mathcal{A} \in \mathbb{R}^{J_1 \times J_2 \times J_1   \times J_2} $ and $\mathcal{B}, \mathcal{C}^T \in \mathbb{R}^{J_1 \times  J_2 \times K_1  \times K_2}.$ As mentioned before, $\mathcal{P}$ and $\mathcal{Q}$ are known as
the reachability and the observability Gramians. As mentioned earlier, since  $\mathcal{P}$ and $\mathcal{Q}$ are weakly
symmetric positive-definite square tensor, then we can obtain the Cholesky-like factors of the two
gramians described as follows
\begin{equation}
	\mathcal{P}=\mathcal{U}*\mathcal{U}^T, \  \  \ \mathcal{Q}=\mathcal{L}*\mathcal{L}^T.  
\end{equation}

Where the tensors $\mathcal{U}$ and $\mathcal{L}$ are of appropriate dimensions, represented in low-rank form. The next step  
involves computing the singular value decomposition (SVD) of $\Phi ( \mathcal{L}^T * \mathcal{U}) = Y \Sigma Z^T $. This decomposition can be represented as follows
%Careful computation is necessary at this stage since tensor $\mathcal{U}$ and
%$\mathcal{L}$ are not of large dimensions. The matrix $ \Sigma $ is a diagonal matrix that contains t. 

\begin{eqnarray}
\Phi (\mathcal{L}^T * \mathcal{U})  & = & Y \Sigma Z^T \nonumber\\
&  =  &  \left[ \begin{array}{l}
	Y_{1} \  \ 	Y_{2} \\
\end{array}\right] \left[ \begin{array}{l}
\Sigma_1  \  \ 	 \\
  \  \  \	\Sigma_2 \\
\end{array}\right] \left[ \begin{array}{c}
	(Z_{1})^T \\
	(Z_{2})^T \\
\end{array}\right] \nonumber\\
	\end{eqnarray}
	
	where $\Sigma_1 \in \mathbb{R}^{r^2 \times r^2}$ and $\Sigma _2 \in \mathbb{R}^{(N^2-r^2) \times (N^2-r^2)}$. As outlined in \cite{gugerc1, mehrman, moore}, a truncation step could
	be established, and by truncating the states that correspond to small Hankel singular values in $\Sigma_2$. Define 
	
	\begin{equation}
		\mathcal{W}_r=\mathcal{L}_1 * \Phi ^{-1}(Y_1)* \Phi^{-1} (\Sigma_1)^{-1/2}, \  \  \  \mathcal{V}_r=\mathcal{U}_1 * \Phi ^{-1}(Z_1)* \Phi^{-1} (\Sigma_1)^{-1/2},
	\end{equation}	
	
	where $\Phi ^{-1}$  is the inverse of the isomorphism defined in Definition \ref{def1}. The matrices $Y_1$ and $Z_1$ are composed of the leading r columns of Y and Z, respectively. We can easily verify that $\mathcal{W}_r ^T* \mathcal{V}_r = \mathcal{I}_r$ and hence that $\mathcal{V}_r* \mathcal{W}_r ^T$ is an oblique projector. The tensor structure of the reduced MLTI system is given as follows
	\begin{equation}
		\mathcal{A}_r=\mathcal{W}_r ^T *(\mathcal{A}*\mathcal{V}_r), \   \   \mathcal{B}_r=\mathcal{W}_r ^T *\mathcal{B},  \  \  \mathcal{C}_r=\mathcal{C}*\mathcal{V}_r.
	\end{equation}

	Regarding the solutions $\mathcal{P}$ and $\mathcal{Q}$ of the
	two Lyapunov equations (\ref{gramians}), we suggest a solution in a factored form via the tensor rational block
	Krylov subspace projection method. This method will be described in the upcoming section.
%In this introduction, we will explore the fundamental concepts behind Tensor Balanced Truncation, its mathematical foundations, and its practical applications. By understanding these elements, we aim to highlight the significance of TBT in modern computational practices and its future potential in various scientific and engineering domains.

\section{Tensor continuous-time Lyapunov equations}
In this section, we discuss the process of obtaining approximate solutions to the tensor continuous  Lyapunov equations of the form

\begin{equation}\label{lyap}
	\mathcal{A}*\mathcal{X}+\mathcal{X}*\mathcal{A}^T+\mathcal{B}*\mathcal{B}^T=\mathcal{O},
\end{equation}

where $\mathcal{A}$ and $\mathcal{B}$ are tensors of appropriate dimensions, and $\mathcal{X}$ is the unknown tensor. Solving this equation is a primary task in the Balanced Truncation model order reduction method for MLTI systems as described in \cite{chen2}. These equations also play a crucial role in control theory; for instance, equation (\ref{lyap}) arises from the discretization of the 2D heat equations with control \cite{chen2, nip}. It is evident that if the dimension of equation (\ref{lyap}) is small, one can transform it into a matrix Lyapunov equation and apply efficient direct techniques as described in \cite{Barr, heyouni, simon}. However, in the case of large-scale tensors, opting for a process based solely on tensors is more beneficial than using matricization techniques, which can be costly in terms of both computation and memory. In the next section, we propose a method to solve continuous Lyapunov tensor equations using tensor Krylov subspace techniques. We use the block Lanczos process based on the tensor rational Krylov subspace, described in the previous section. 

\subsection{Tensor rational block Lanczos method for continuous-Lyapunov equation} 
	
We recall the tensor Lyapunov equation that we are interested in

\begin{equation}\label{lyap2}
	\mathcal{A}*\mathcal{X}+\mathcal{X}*\mathcal{A}^T+\mathcal{B}*\mathcal{B}^T=\mathcal{O},
\end{equation}

where $\mathcal{A} \in \mathbb{R}^{J_1 \times  J_2 \times J_1  \times J_2}, \mathcal{B} \in \mathbb{R}^{J_1 \times  J_2 \times K_1 \times K_2}$ and $\mathcal{X} \in \mathbb{R}^{J_1  \times J_2 \times J_1   \times J_2}$ is the unknown tensor. Analogously to the matrix continuous-Lyapunov equation \cite{horn}, we state that (\ref{lyap2}) has a unique solution if the
following condition is satisfied

$$\lambda _{j_1 j_2 j_1 j_2}(\mathcal{A})+\bar{\lambda} _{j_1 j_2 j_1 j_2}(\mathcal{A}) \neq 0, \  \  j_1=1,\ldots,J_1; j_2=1,\ldots,J_2$$
 where $\lambda _{j_1 j_2 j_1 j_2}(\mathcal{A})$ and it’s conjugate $\bar{\lambda} _{j_1 j_2 j_1 j_2}(\mathcal{A})$ are eigenvalues of $\mathcal{A}$ defined in Definition \ref{eigenv}.
 
 Next, we describe the procedure of constructing the approximation by using the tensor rational block Lanczos process defined in Algorithm \ref{alg2}. We seek an approximate solution $\mathcal{X}_m$ to $\mathcal{X}$ that satisfies the
 continuous Lyapunov equation (\ref{lyap2}). This approximate solution is given by
 \begin{equation}\label{appr}
\mathcal{X}_m = \mathcal{V}_m* \mathcal{Y}_m * \mathcal{V}_m^T,
 \end{equation}

The tensor $\mathcal{Y}_m \in \mathbb{R}^{K_1 \times  m K_2 \times K_1  \times m K_2}$ such that the following Galerkin condition must be satisfied
\begin{equation}\label{galerkin}
	\mathcal{W}_m ^T* (\mathcal{R}_m* \mathcal{W}_m)=\mathcal{O},
\end{equation}

where $\mathcal{V}_m$ and $\mathcal{W}_m$ are the tensors obtained after running Algorithm \ref{alg2} to the triplet $(\mathcal{A}, \mathcal{B}, \mathcal{C}^T)$. $\mathcal{R}_m$ is the residual tensor given by $\mathcal{R}_m =\mathcal{A}*\mathcal{X}_m+\mathcal{X}_m*\mathcal{A}^T+\mathcal{B}*\mathcal{B}^T.$ Using the bi-orthogonality condition (i.e., $\mathcal{W}_m ^T * \mathcal{V}_m = \mathcal{I}_m$) and developing the condition (\ref{galerkin}),
we find that the tensor $\mathcal{Y}_m$ satisfying the following low dimensional Lyapunov tensor equation
\begin{equation}\label{LRE}
\mathcal{T}_m * \mathcal{Y}_m+\mathcal{Y}_m*\mathcal{T}_m^T+\mathcal{B}_m*\mathcal{B}_m^T=\mathcal{O},
\end{equation}
where $\mathcal{T}_m= \mathcal{W}_m ^T * \mathcal{A} * \mathcal{V}_m $ and $\mathcal{B}_m = \mathcal{W}_m ^T * \mathcal{B}.$

The following result shows an efficient way to compute the residual error in an
efficient way.

\begin{theorem}\label{res}
Assume that Algorithm \ref{alg2} has been executed for m iterations. Then	

\begin{equation}
	\Vert \mathcal{R}_m \Vert \leq 2 \Vert \Gamma_{m, A} \mathcal{X}_m \mathcal{V}_m ^T \Vert
\end{equation}

where $\Vert . \Vert $ denotes the Frobenius tensor norm, as defined in equation  (\ref{def3}) and $\Gamma_{m, A}$ is defined in equations \ref{condit2}.
	
	\end{theorem}
	
	\begin{proof}
		
By using the definition of the approximation $\mathcal{X}_m$ given in equation (\ref{appr}) and the decompositions  provided in equation (\ref{condit2}) from Algorithm  \ref{alg2}, we obtain		

\begin{eqnarray*}
	\mathcal{R}_m & =  & \mathcal{A}*\mathcal{X}_m+\mathcal{X}_m*\mathcal{A}^T+\mathcal{B}*\mathcal{B}^T  \\ \nonumber
		&  =  &  \mathcal{A} * (\mathcal{V}_m * \mathcal{Y}_m * \mathcal{V}_m ^T)+ (\mathcal{V}_m * \mathcal{Y}_m * \mathcal{V}_m ^T)*\mathcal{A}^T + \mathcal{B}*\mathcal{B}^T \\ \nonumber
		& =  &  (\mathcal{V}_{m}* \mathcal{T}_m + \Gamma _{m, A}) * \mathcal{Y}_m *\mathcal{V}_m ^T + \mathcal{V}_{m}* \mathcal{Y}_{m}*(\mathcal{T}_{m}^T* \mathcal{V}_m^T + \Gamma _{m, A} ^T) + \mathcal{V}_m * (\mathcal{B}_m * \mathcal{B}_m^T)*\mathcal{V}_m^T \\ \nonumber
		&  = & \mathcal{V}_m *  \left[ \mathcal{T}_m * \mathcal{Y}_m+\mathcal{Y}_m*\mathcal{T}_m^T+\mathcal{B}_m*\mathcal{B}_m^T   \right]  * \mathcal{V}_m^T +  \Gamma _{m, A} * \mathcal{Y}_m * \mathcal{V}_m ^T + \mathcal{V}_m * \mathcal{Y}_m * \Gamma _{m, A}^T \\ \nonumber
		& = & \Gamma _{m, A} * \mathcal{Y}_m * \mathcal{V}_m ^T + \mathcal{V}_m * \mathcal{Y}_m * \Gamma _{m, A}^T, \\ \nonumber
\end{eqnarray*}
where the last equation is obtained using the fact that $ \mathcal{V}_m * \mathcal{W}_m^T * \mathcal{B} = \mathcal{B}$ and equation (\ref{LRE}).

As $\mathcal{Y}_m$ is a symmetric matrix, it follows that

$$	\Vert \mathcal{R}_m \Vert \leq 2 \Vert \Gamma_{m, A} \mathcal{X}_m \mathcal{V}_m ^T \Vert .$$

	\end{proof}		
	
\begin{remark}
For an efficient computation, the approximation $\mathcal{X}_m$ can be expressed in a factored form. Then we begin by applying  Singular Value Decomposition (SVD) to $\mathcal{Y}_m$, i.e., $\mathcal{Y}_m = \mathcal{P} * \Sigma * \mathcal{Q}^T$. Next, we consider a tolerance {\tt dtol} and define $\mathcal{P}_r$ and $\mathcal{Q}_r$ as the first r columns of $\mathcal{P}$ and $\mathcal{Q}$, respectively, corresponding to the first r singular values whose magnitudes exceed {\tt dtol}. By setting $\Sigma_r =diag(\sigma_1,\ldots, \sigma_r)$, we approximate $\mathcal{Y}_m$ as $\mathcal{Y}_m \approx \mathcal{P}_r * \Sigma_r * \mathcal{Q}_r ^T$. The factorization then proceeds as follows:

\begin{equation}
	\mathcal{X}_m \approx \mathcal{V}_m *( \mathcal{P}_r * \Sigma_r * \mathcal{Q}_r ^T )*\mathcal{V}_m ^T \approx \mathcal{Z}_1*\mathcal{Z}_2 ^T,
\end{equation}
where $\mathcal{Z}_1 = \mathcal{V}_m * (\mathcal{P}_r * (\Sigma_r) ^{1/2})$ and $\mathcal{Z}_2^T=  ((\Sigma_r) ^{1/2}*\mathcal{Q}_r ^T) * \mathcal{V}_m ^T.$
\end{remark}	

The following algorithm describe all the results given in this section.

\begin{algorithm}[h]
	\caption{The Lyapunov tensor rational block Lanczos  algorithm {\tt(LTRBLA)} }
	\label{alg3}
	\begin{enumerate}
		\item \textbf{Input:} $\mathcal{A}\in\mathbb{R}^{J_1 \times J_2 \times J_1 \times J_2},\ \mathcal{B}, \mathcal{C}^T\in\mathbb{R}^{J_1 \times J_2 \times K_1 \times K_2}$, tolerances $\epsilon,$ {\tt dtol} and a fixed integer $m_{max}$  of maximum iterations.
		\item \textbf{Output:}  The approximate solution $\mathcal{X}_m = \mathcal{Z}_1 * \mathcal{Z}_2 ^T$.
		\item \textbf{For} $j=m,\ldots,m_{max}$ \textbf{do}
		\item \ \  Use Algorithm \ref{alg2} to built $\mathcal{V}_{m}$ and $\mathcal{W}_{m}$ bi-orthonormal tensors associated to the tensor
		rational Krylov subspaces $\mathcal{K}_m (\mathcal{A}, \mathcal{B}, \Sigma_m)$ and $\mathcal{K}_m (\mathcal{A}^T, \mathcal{C}^T, \Sigma_m)$ defined in ( \ref{ratkrylov}) and compute the tensor $\mathcal{T}_m$ using decompositions (\ref{condit2}).
		\item \ \  Solve the low-dimensional continuous Lyapunov equation (\ref{LRE}) using the MATLAB
		function {\tt dlyap}. 
		\item \ \  Compute the residual norm $\Vert \mathcal{R}_m \Vert $ using Theorem \ref{res}, and if it is less than $\epsilon$, then
		 
		 \ \ \ \ \  a) compute the SVD of $\mathcal{Y}_m = \mathcal{P} * \Sigma * \mathcal{Q}^T$ where $\Sigma=diag(\sigma_1,\ldots,\sigma_m)$.
		
		\ \ \ \ \  b) determine r such that $\sigma_{r+1} < {\tt dtol} \leq \sigma_{r}$. Set $\Sigma _r=diag(\sigma_1,\ldots,\sigma_r)$ and compute
		 $$\mathcal{Z}_1 = \mathcal{V}_m * (\mathcal{P}_r * (\Sigma_r) ^{1/2}) \  \ \ and \ \ \  \mathcal{Z}_2=  \mathcal{V}_m * (\mathcal{Q}_r * (\Sigma_r) ^{1/2}).$$
		\item  \textbf{endFor}
		\item  \textbf{return  $\mathcal{X}_m$}
	\end{enumerate}
\end{algorithm}

\subsection{The coupled Lyapunov equations}
 Now, let us see how to solve the coupled Lyapunov equation (\ref{gramians}) by using the tensor rational block Lanczos process.

 First, we need to solve the  following two low-dimensional Lyapunov equations

\begin{equation}\label{2LRE}
	\left\{ 
	\begin{array}{c c c}
		\mathcal{T}_m * \mathcal{X}_{m}+\mathcal{X}_{m}*\mathcal{T}_m^T+\mathcal{B}_m*\mathcal{B}_m^T=\mathcal{O}\\ 	
		\mathcal{T}_m ^T * \mathcal{Y}_{m}+\mathcal{Y}_{m}*\mathcal{T}_m+\mathcal{C}_m ^T*\mathcal{C}_m=\mathcal{O}.\\
	\end{array}
	\right. 
	%\end{array}%
	%\]
\end{equation}
 We then form the approximate solutions $\mathcal{P}_m = \mathcal{V}_m* \mathcal{X}_{m} * 	\mathcal{V}_m^T,$ and $\mathcal{Q}_m = \mathcal{W}_m* \mathcal{Y}_{m} * 	\mathcal{W}_m^T,$  where $\mathcal{V}_m$ and $\mathcal{W}_m$ are the tensors obtained after running Algorithm \ref{alg2} on the triplet $(\mathcal{A}, \mathcal{B}, \mathcal{C}^T)$ and $\mathcal{T}_m=\mathcal{W}_m^T * \mathcal{A}*\mathcal{V}_m$. We summarize the coupled Lyapunov tensor rational block Lanczos algorithm as follows.
 
 \begin{algorithm}[h]
 	\caption{The Coupled Lyapunov Tensor Rational Block Lanczos  algorithm {\tt(CLTRBLA)}  }
 	\label{alg4}
 	\begin{enumerate}
 		\item \textbf{Input:} $\mathcal{A}\in\mathbb{R}^{J_1 \times J_2 \times J_1 \times J_2},\ \mathcal{B}, \mathcal{C}^T\in\mathbb{R}^{J_1 \times J_2 \times K_1 \times K_2}$.
 		\item \textbf{Output:}  The approximate solutions $\mathcal{P}_m$ and $\mathcal{Q}_m$.
 		\item Apply Algorithm 3 to the triplet $(\mathcal{A}, \mathcal{B}, \mathcal{C}^T)$.
 		\item The approximate solutions are represented as the tensor products:
$$\mathcal{P}_m = \mathcal{V}_m* \mathcal{X}_{m} * 	\mathcal{V}_m^T,  \  \  \textit{and} \   \   \    \mathcal{Q}_m =\mathcal{W}_m* \mathcal{Y}_{m} * 	\mathcal{W}_m^T.$$
 	\end{enumerate}
 \end{algorithm}

\section{Numerical experiments}

In this section, we give some experimental results to show
the effectiveness of the tensor rational block Arnoldi process ({\tt TRBA}) and the tensor rational block Lanczos process ({\tt TRBL}),  when applied to reduce the
order of multidimensional linear time invariant systems. The numerical results were obtained using MATLAB R2016a on a computer with Intel core i7 at 2.6GHz and 16 Gb of RAM. We need to mention that all the algorithms described here
have been implemented based on the MATLAB tensor toolbox developed  by Kolda et al., \cite{kolda}.

{\bf Example 1.}  In this example, we applied the {\tt TRBA} and the {\tt TRBL} processes to reduce the order of multidimensional
linear time invariant (MLTI) systems.

{\bf Example 1.1.} For the first experiment, we used the following data:

$\bullet  \mathcal{A} \in \mathbb{R}^{N \times N \times N \times N}$, with N=100. Here, A is constructed from a tensorization of a triangular
matrix $ A \in \mathbb{R}^{N N \times N N}$ constructed using the MATLAB function {\tt spdiags}.

$\bullet  \mathcal{B}, \mathcal{C}^T  \in \mathbb{R}^{N \times N \times K_1 \times K_2}$ are chosen as sparse and random tensors, respectively, with $K_1 =3,  K_2 = 4$.

{\bf Example 1.1.} For the second experiment, we consider  the evolution of heat distribution in a
solid medium over time. The partial differential equation that describes this evolution is known as the 2D heat equation, given by the following equation

\begin{equation*}
	\left\{ 
	\begin{array}{c c c}
	\dfrac{	\partial}{\partial t} \phi(t,x)  & = & c^2 	\dfrac{	\partial ^2}{\partial x ^2} \phi(t,x)+\delta(x) u_t ,(\textit{square} D=[-\pi ^2, \pi]) \  \    \    x \in D\\ 
		\phi(t,x)  & = & 0  \   \   \  x \in \partial D,\\
	\end{array}
	\right. 
	%\end{array}%
	%\] 
\end{equation*}

$\bullet$  The tensor $\mathcal{A} \in \mathbb{R}^{N \times N \times N \times N}$, with N=80 is the tensorization of $\dfrac{c^2 \delta t}{h^2} \Delta_{dd} \in \mathbb{R}^{N^2 \times N^2} $, where $\Delta _{dd}$  d is the discrete Laplacian on a rectangular grid with  Dirichlet boundary condition
(see \cite{chen2} for more details).

$\bullet  \mathcal{B}, \mathcal{C}^T  \in \mathbb{R}^{N \times N \times K_1 \times K_2}$ are chosen to be random tensors  with $K_1 =3,  K_2 = 4$.

\begin{figure}[!h]
	\begin{center}
		\includegraphics[height=3in ,width=2.75in]{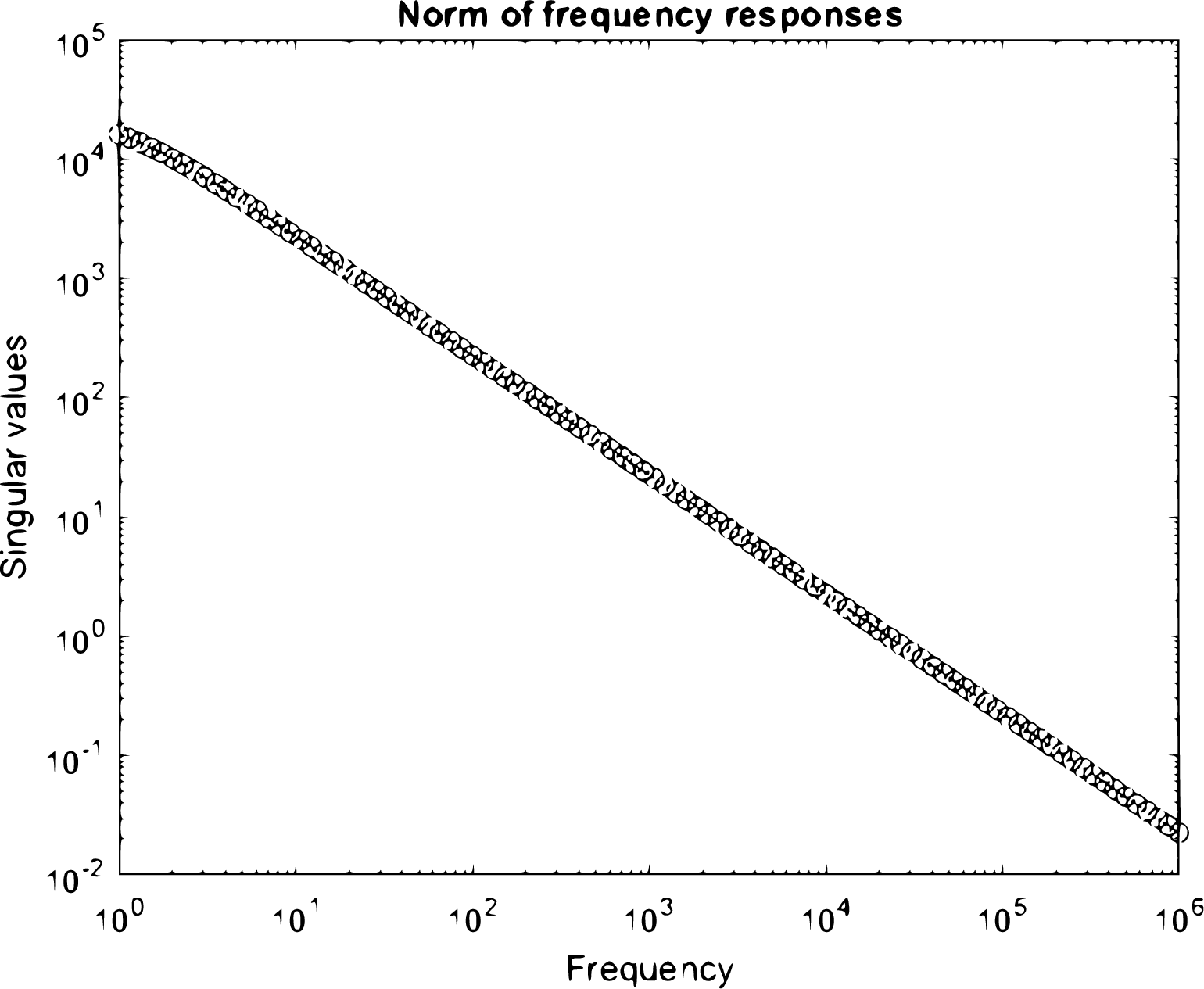}
		\includegraphics[height=3in ,width=2.75in]{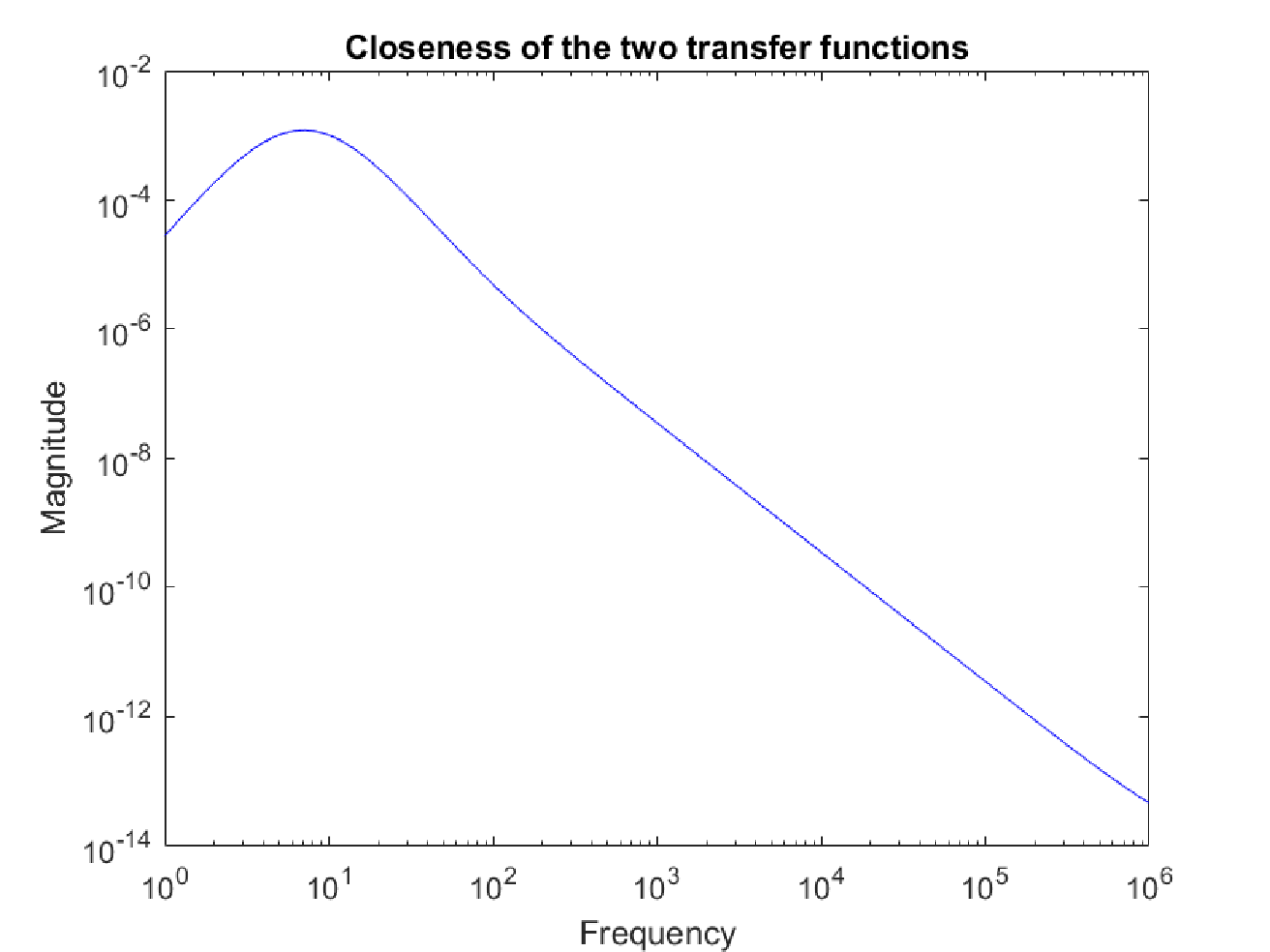}
	\end{center}
	\caption{{\bf Example 1.1.} Left: $\Vert {\mathcal{F}(j \omega)} \Vert _{2}$ and it's approximation $\Vert \mathcal{F}_m(j \omega) \Vert _{2}$. Right: the exact error $\Vert \mathcal{F}(j  \omega)- \mathcal{F}_m(j \omega) \Vert _{2}$ using the {\tt TRBA} algorithm  with $m=5$.}\label{fig1.1}
\end{figure}

\begin{figure}[!h]
	\begin{center}
		\includegraphics[height=3in ,width=2.75in]{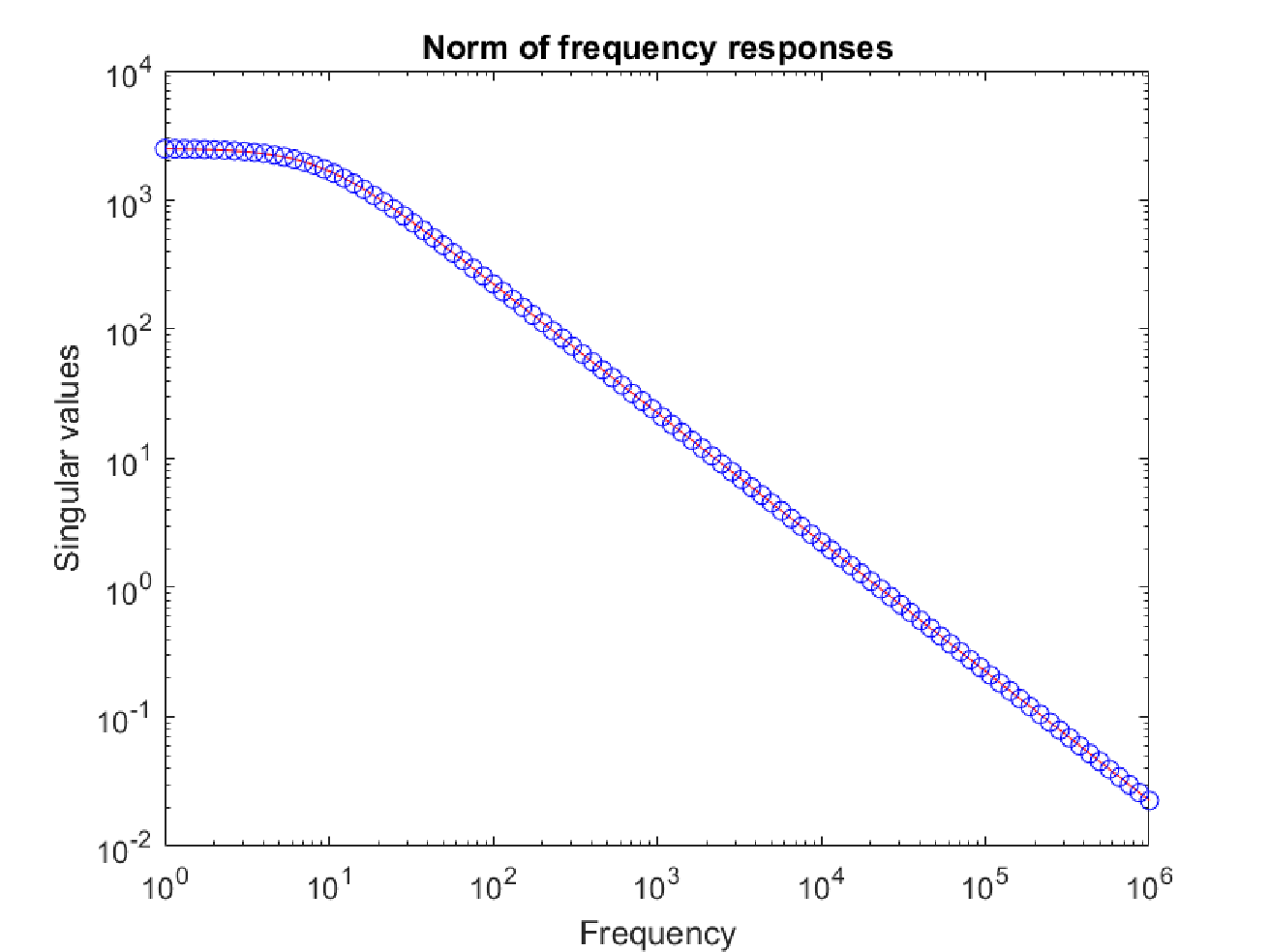}
		\includegraphics[height=3in ,width=2.75in]{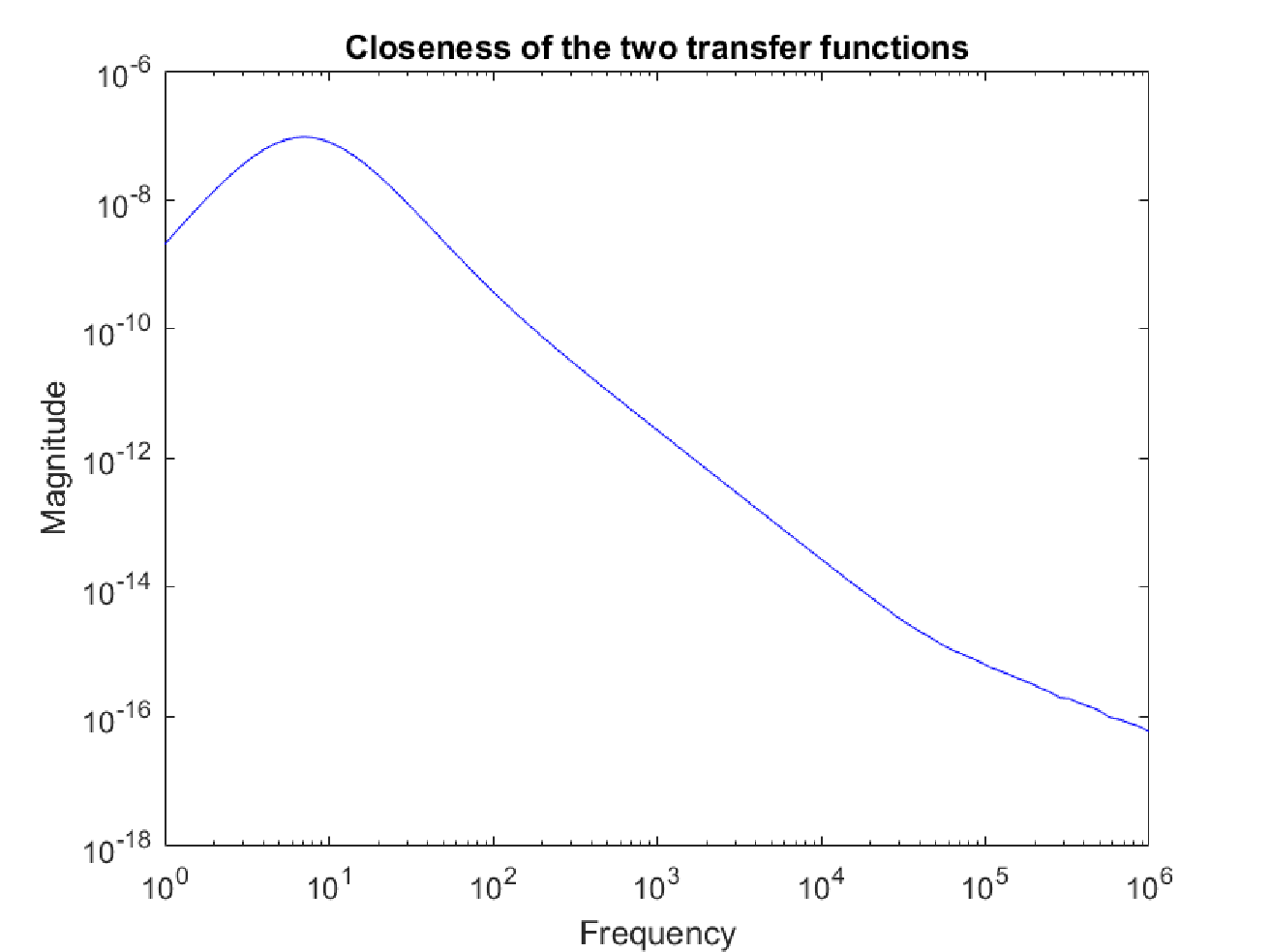}
	\end{center}
	\caption{{\bf Example 1.1.} Left: $\Vert {H(j \omega)} \Vert _{2}$ and it's approximation $\Vert H_m(j \omega) \Vert _{2}$. Right: the exact error $\Vert H(j  \omega)- H_m(j \omega) \Vert _{2}$ using the {\tt TRBL} algorithm  with $m=10$.}\label{fig1.2}
\end{figure}

\begin{figure}[!h]
	\begin{center}
		\includegraphics[height=3in ,width=2.75in]{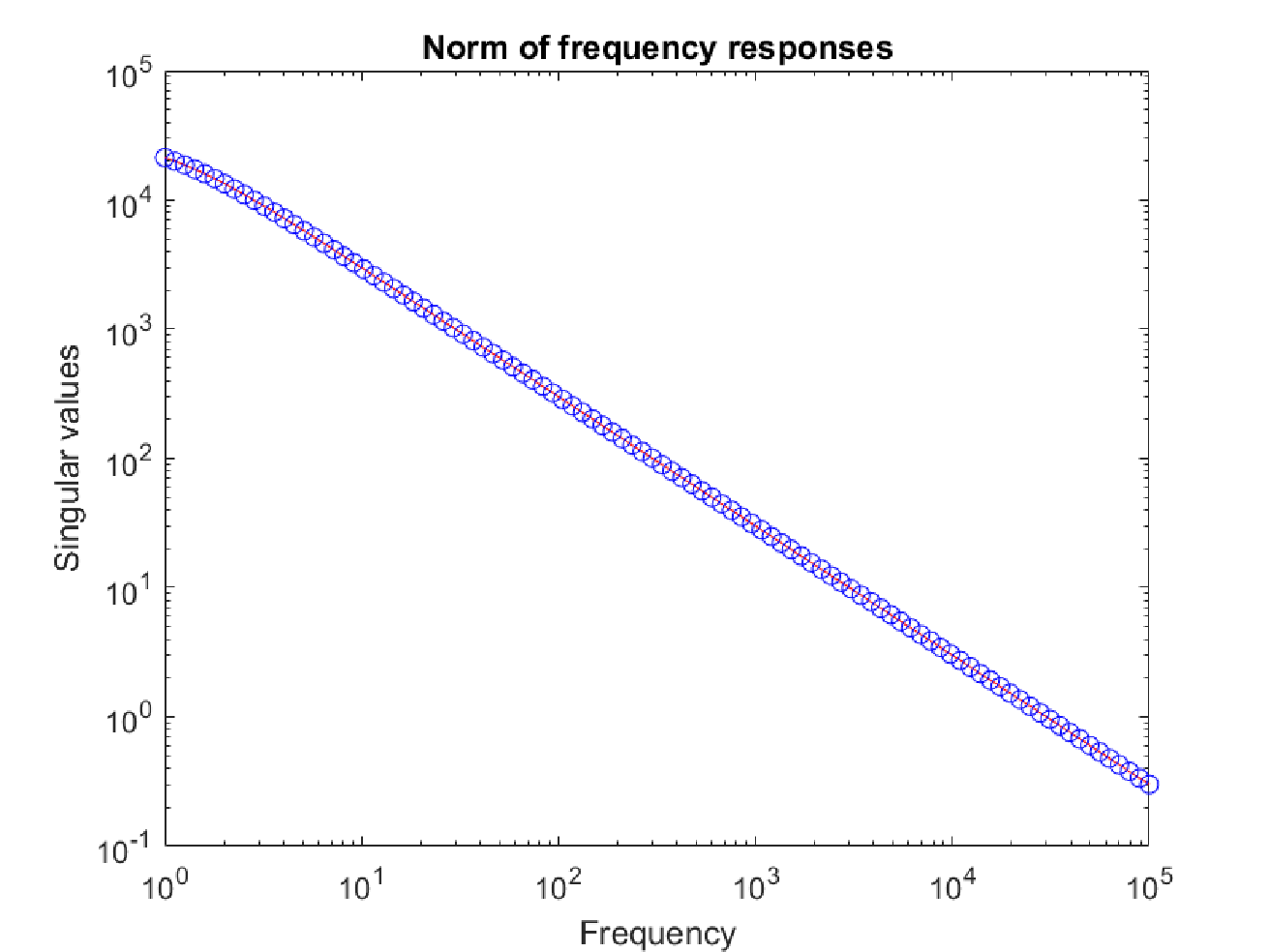}
		\includegraphics[height=3in ,width=2.75in]{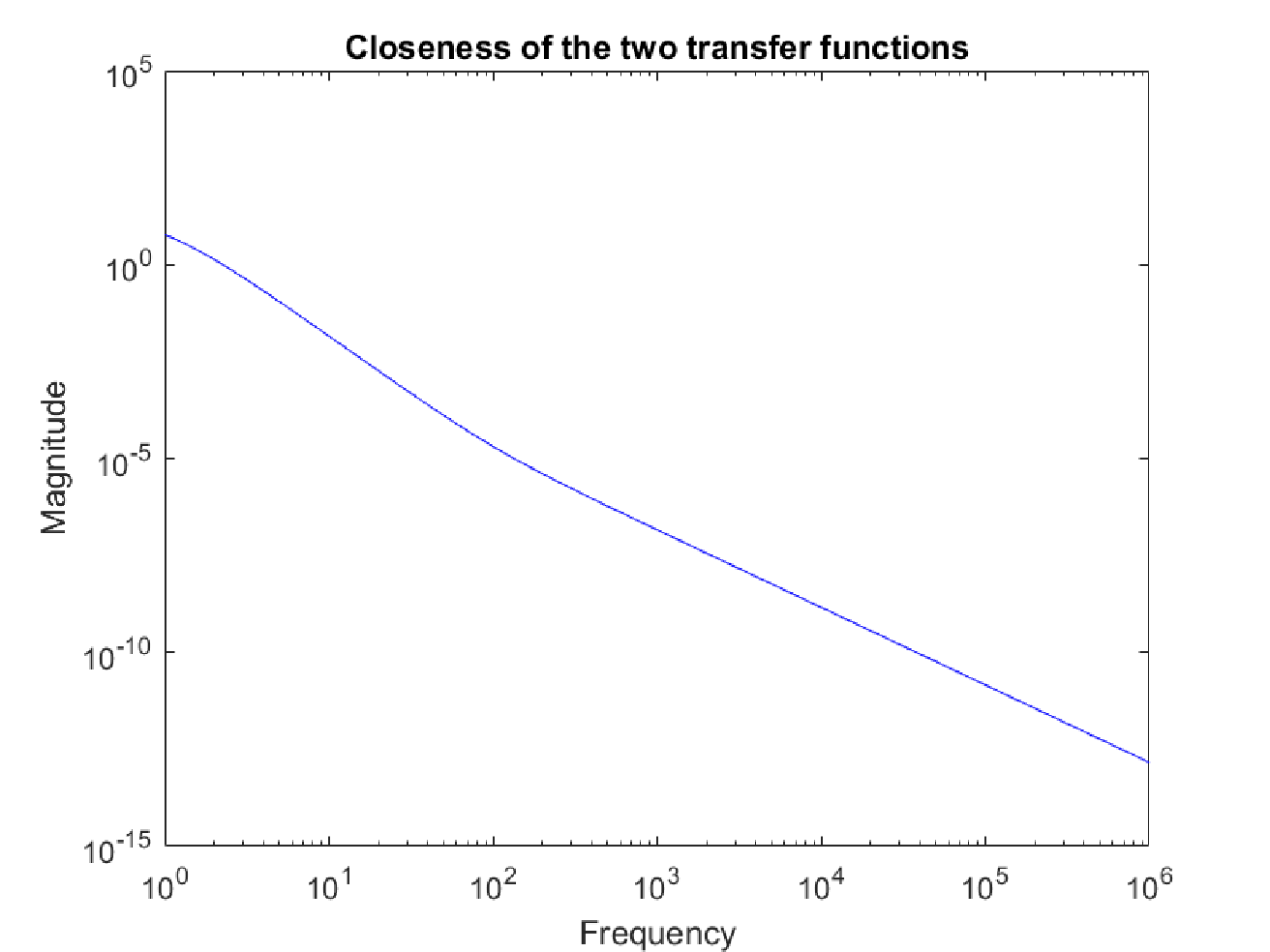}
	\end{center}
	\caption{{\bf Example 1.2.} Left: $\Vert {H(j \omega)} \Vert _{2}$ and it's approximation $\Vert H_m(j \omega) \Vert _{2}$. Right: the exact error $\Vert H(j  \omega)- H_m(j \omega) \Vert _{2}$ using the {\tt TRBA} algorithm  with $m=10$.}\label{fig2.1}
	\end{figure}

\begin{figure}[!h]
		\begin{center}
			\includegraphics[height=3in ,width=2.75in]{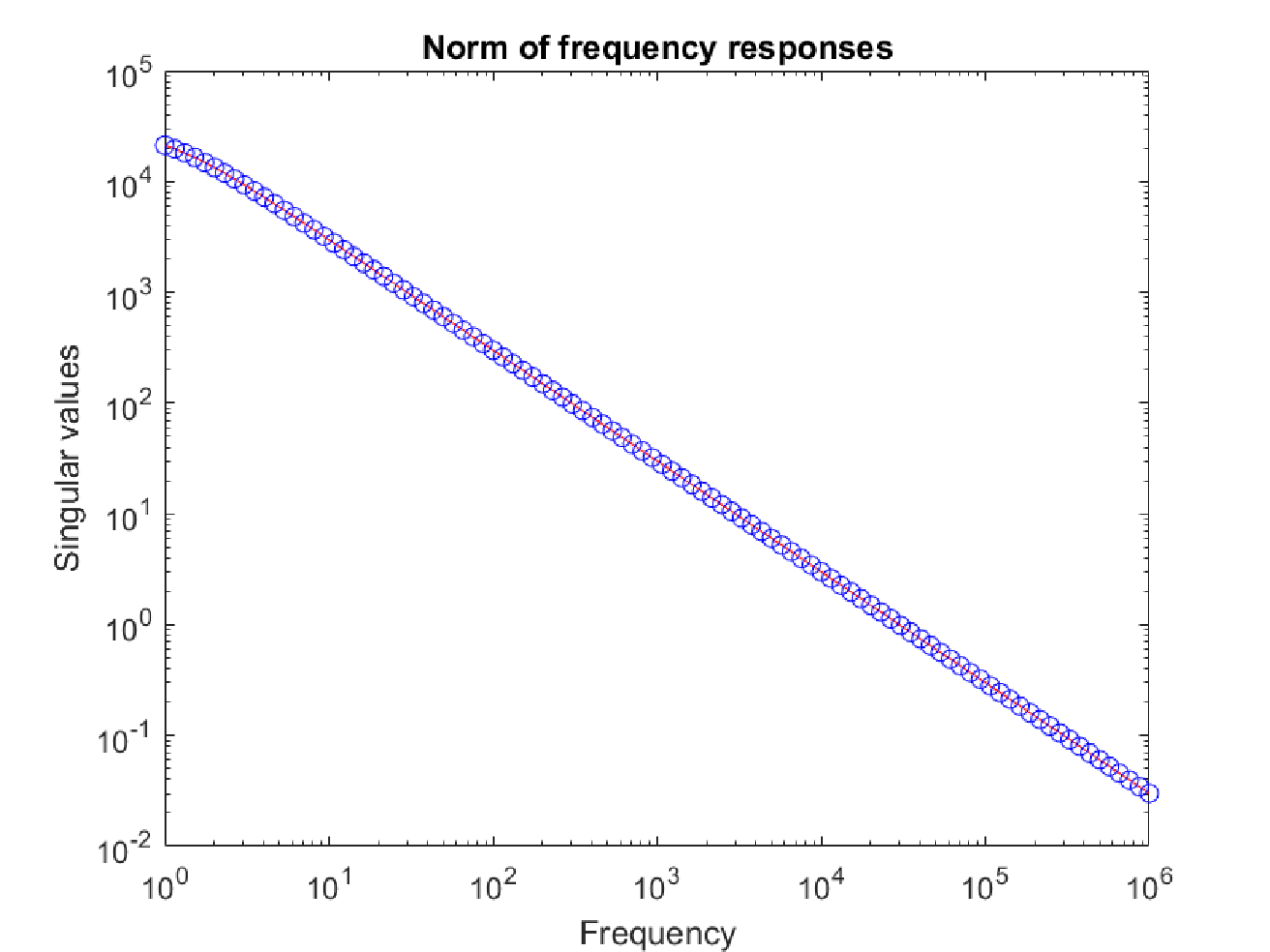}
			\includegraphics[height=3in ,width=2.75in]{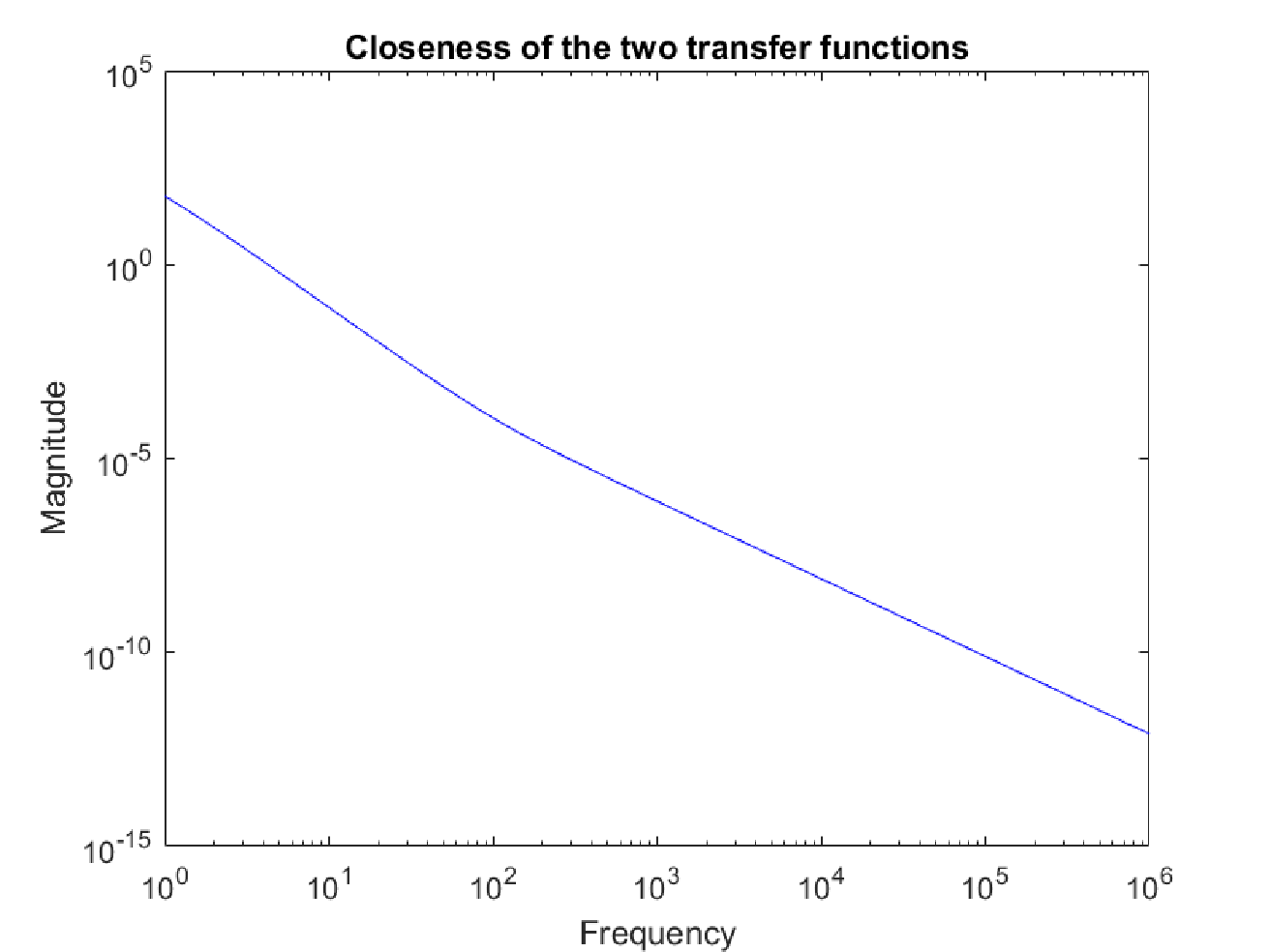}
		\end{center}
		\caption{{\bf Example 1.2.} Left: $\Vert {\mathcal{F}(j \omega)} \Vert _{2}$ and it's approximation $\Vert \mathcal{F}_m(j \omega) \Vert _{2}$. Right: the exact error $\Vert \mathcal{F}(j  \omega)- \mathcal{F}_m(j \omega) \Vert _{2}$ using the {\tt TRBL} algorithm  with $m=15$.}\label{fig2.2}
		\end{figure}

 The performance of the methods is shown in the plots of Figs. \ref{fig1.1}–\ref{fig2.2}, which report the frequency response of the original system (circles plot) compared to the frequency response of its approximation (solid plot). The right plots
 of these figures represent the exact error  $\Vert \mathcal{F}(j  \omega)- \mathcal{F}_m(j \omega) \Vert _{2}$ for different frequencies and different  small space dimension.
 
%which report the Bode plot of the transfer function (left plots) and
%the associated errors (right plots) for the reference (small) space

\newpage

{\bf Example 2.} For this experiment, we show the results obtained from the tensor Balanced
Truncation method explained in Section 5, to reduce the order of {\tt MLTI} systems. The following data were used:

$ \bullet $ The tensor $\mathcal{A} \in \mathbb{R}^{N \times N \times N \times N}$, with N = 100, is constructed from a tensorization of a triangular
matrix $ A \in \mathbb{R}^{N N \times N N}$ constructed using the MATLAB function {\tt spdiags}.

$\bullet  \mathcal{B}, \mathcal{C}^T  \in \mathbb{R}^{N \times N \times K_1 \times K_2}$ are chosen to be random tensors  with $K_1 = K_2 = 3$.

$ \bullet $ We use Algorithm 3 to solve the two continuous Lyapunov equations (\ref{gramians}) by setting \\
m = 30 (i.e., the maximum number of iterations) and the tolerance $\epsilon = 10^{-8}$.

\begin{figure}[!h]
	\begin{center}
		\includegraphics[height=3in ,width=2.75in]{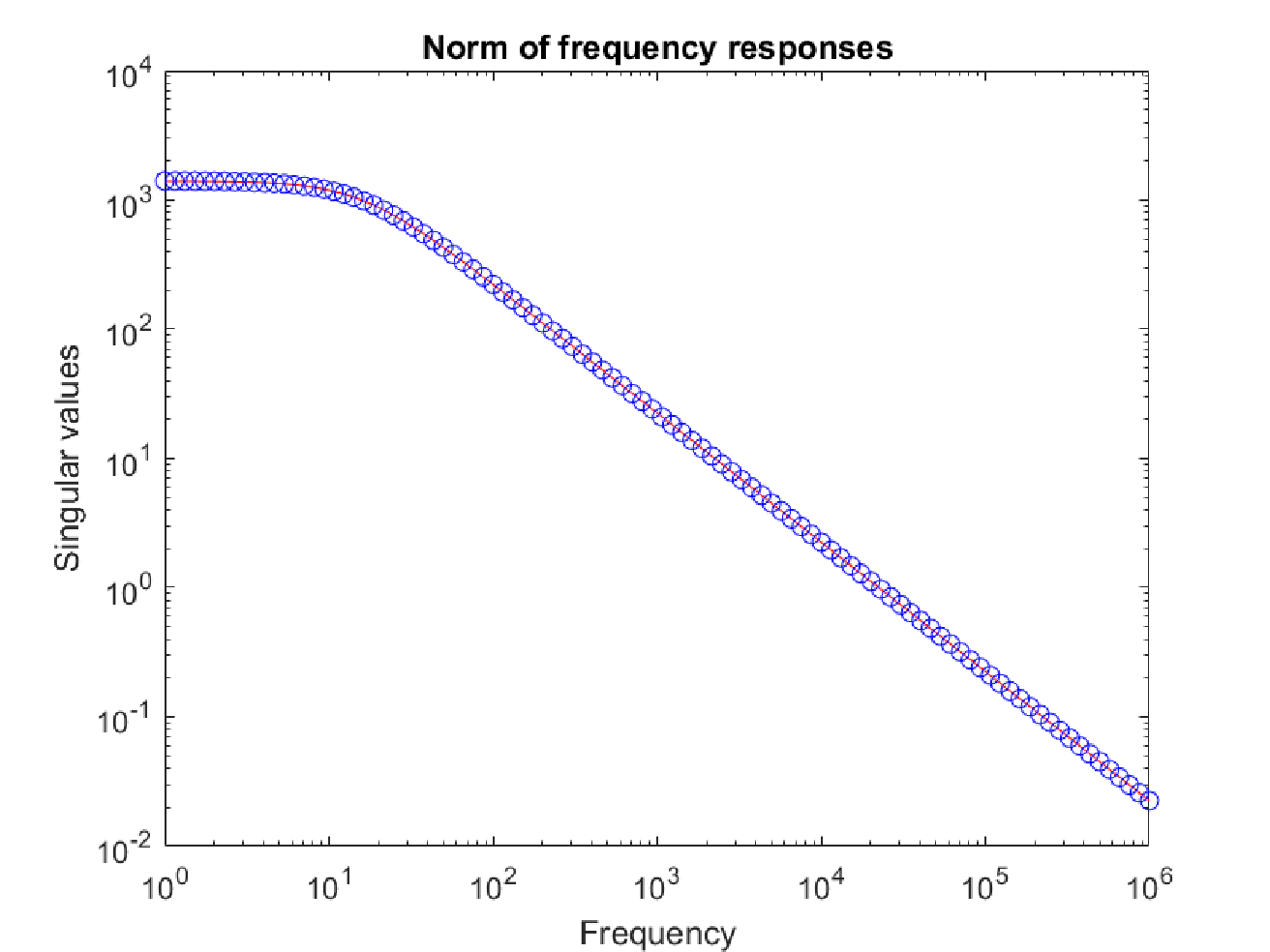}
		\includegraphics[height=3in ,width=2.75in]{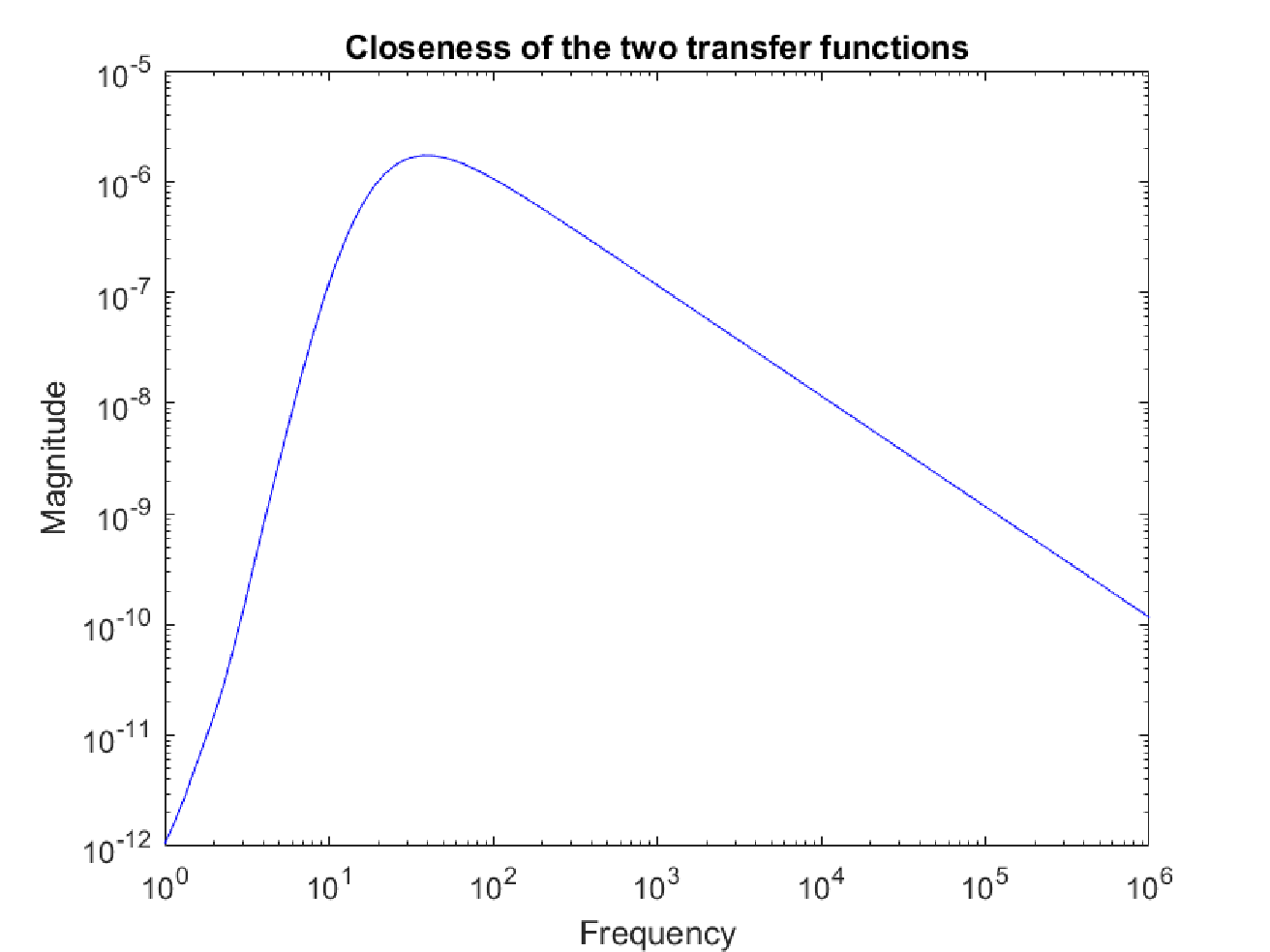}
	\end{center}
	\caption{Left: $\Vert {\mathcal{F}(j \omega)} \Vert _{2}$ and it's approximation $\Vert \mathcal{F}_m(j \omega) \Vert _{2}$. Right: the exact error $\Vert \mathcal{F}(j  \omega)- \mathcal{F}_m(j \omega) \Vert _{2}$ for m=5.}\label{fig3}
\end{figure}

The frequency response (solid plot) of this test example  is given in the left of Fig. 7.5 and is compared
to the frequency responses of the 5 order approximation (circles plot). The exact error $\Vert \mathcal{F}(j  \omega)- \mathcal{F}_m(j \omega) \Vert _{2}$ produced by this  process is shown in the right of Fig. \ref{fig3}.

{\bf Example 3.} For the last experiment, we compared the performance of the tensor classical block Lanczos ({\tt TCBL})  and the tensor rational block Lanczos ({\tt TRBL}) algorithms when applied to solve high order continuous-time
Lyapunov tensor equations.  The following data were used:

$ \bullet $ The tensor $\mathcal{A} \in \mathbb{R}^{N \times N \times N \times N}$ is constructed from a tensorization of a triangular
matrix $ A \in \mathbb{R}^{N N \times N N}$ constructed using the MATLAB function {\tt spdiags}.

$\bullet  \mathcal{B}, \mathcal{C}^T  \in \mathbb{R}^{N \times N \times K_1 \times K_2}$ are chosen to be random tensors.

$ \bullet $ We use Algorithm 3 to solve the two continuous Lyapunov equations (\ref{gramians}) by setting 
m = 30 (i.e., the maximum number of iterations) and the tolerance $\epsilon = 10^{-8}$.

 As observed from Table \ref{tab1}, {\tt TRBL} process is more effective than the {\tt TCBL} process, and the number of iterations required for convergence is small. However, the {\tt TCBL} algorithm  performs better in terms of computation time.

\begin{table}[!h]
	\begin{center}
		\caption{Results obtained from the two methods for solving high order coupled Lyapunov tensor equations.}\label{tab1}
		\begin{tabular}{l|c|c|c|c}
			\hline
			$(N,K_1,K_2)$ & \;\;  Algorithm & \;\;   Iter  & \;\;  Res & \;\; time (s) \\  
			\hline
			{\tt (80,3,3)} & {\tt TCBL} & 11 & $4.02 \times 10^{-9}$ &  15.42 \\ 
			                     & {\tt TRBL} &  6 & $9.71 \times 10^{-10}$ &  138.24\\
			\hline
			{\tt (80,3,4)} & {\tt TCBL} & 11 & $4.52 \times 10^{-9}$ &  18.80 \\ 
			                     & {\tt TRBL} &  6 & $8.28 \times 10^{-10}$ &  142.5\\
			\hline
                           {\tt (100,3,3)} & {\tt TCBL} & 12 & $3.96 \times 10^{-9}$ &  40 \\
                                                 & {\tt TRBL} &  7 & $8.07 \times 10^{-9}$ &  206.14\\
                            \hline
                           {\tt (100,3,4)} & {\tt TCBL} & 11 & $2.23 \times 10^{-8}$ &  47.85 \\
                                                 & {\tt TRBL} &  7& $4.28 \times 10^{-9}$ &  403.37\\
		
		\end{tabular}
	\end{center}
\end{table}

\section{Conclusion}
In this paper, we introduce the tensor rational block Krylov subspace methods based on Arnoldi and Lanczos algorithms.  We extended the application of these methods to reduce the order of multidimensional linear time-invariant (MLTI) systems and to solve  continuous-time Lyapunov tensor equations. Moreover, we derived some theoretical results concerning the error estimations between the original and the reduced transfer functions  and the residual of continuous-time Lyapunov tensor equations. We
presented an adaptive approach for choosing the interpolation points  used to construct the tensor rational Krylov subspaces.  Finally, numerical results are given to confirm the effectiveness of the proposed methods.

\end{document}